\documentclass{siamart171218}
\usepackage[english]{babel}
\usepackage[margin=3cm]{geometry}
\usepackage{algpseudocode}
\usepackage[]{subfig}
\usepackage{amsmath,amssymb,mathtools}
\usepackage{cite,enumerate,float,pgfplots}
\usepackage{courier,dsfont}
\usepackage[utf8]{inputenc}
\usepackage{enumitem}
\newsiamremark{remark}{Remark}

\usepackage{xcolor}

\title{Adjoint-based Adaptive Multi-Level Monte Carlo for Differential Equations}
\author{Jehanzeb Chaudhry\thanks{The University of New Mexico, Albuquerque, NM 87131, USA
(\email{jehanzeb@unm.edu},\email{zstevens@unm.edu}).}
\and Zachary Stevens\footnotemark[1]}
\date{\today}

\begin{document}

\maketitle

\begin{abstract}
    We present a multi-level Monte Carlo (MLMC) algorithm with adaptively refined meshes and accurately computed stopping-criteria utilizing adjoint-based \textit{a posteriori} error analysis for differential equations. This is in contrast to classical MLMC algorithms that use either a hierarchy of uniform meshes or adaptively refined meshes based on  Richardson extrapolation, and employ a stopping criteria that relies on assumptions on the convergence rate of the MLMC levels. This work develops two adaptive refinement strategies for the MLMC algorithm. These strategies are based on a decomposition of  an error estimate of the MLMC bias and utilize variational analysis, adjoint problems and computable residuals.
\end{abstract}


\section{Introduction}

Uncertainty quantification (UQ) plays a vital role in predictive computational science and numerical analysis \cite{UQ_SC,Oden2017}. Many physical phenomena are modeled by differential equations, often with uncertain inputs such as parameter values, boundary/initial conditions or the geometry of the domain. UQ is used to characterize, quantify, or control the uncertainty in the outputs of the model, based on information about the uncertainty of the input data \cite{Smith2013,xiu2010numerical}. The uncertainty of the output is quantified by approximating the probability distribution or some statistic (e.g. expectation, standard deviation, mean squared error) of the random variable corresponding to the solution of the model. The UQ of many complex applications presents a challenge due to high-dimensional parameter spaces, non-smooth model responses, and non-linearities.

Monte Carlo (MC) methods form a broad approach for UQ and have been used as a classical tool in predictive computational sciences \cite{Heinrich_MLMC,Fishman2013,McClarren18}. MC methods are simple to implement and are used across a wide range of disciplines such as high-energy physics, finance, particle transport, inter-particle collision, statistical mechanics and numerous other engineering applications \cite{NMSU_MCold,AMJ_MC}. However, obtaining an accurate estimate is often computationally expensive in the context of partial differential equation (PDE) models \cite{giles2015}. Several adaptations and expansions have been made to the MC method to improve its cost effectiveness including the Quasi-Monte Carlo (QMC), Markov chain Monte Carlo (MCMC), multi-fidelity Monte Carlo (MFMC), and multi-level Monte Carlo (MLMC) methods \cite{caflisch_1998,Cliffe2011,SchStTec_QMC,DKPS_QMC,AH_Markov,Giles2008}. Adaptive MLMC methods \cite{KornYou_2018}, which employ adaptive refinement strategies for pathwise problems, have also been investigated.

Many widely used MLMC algorithms rely on a number of assumptions \cite{giles2015,Heinrich_MLMC}, hence the reliability of these methods in predictive science requires quantifying the error in the output of these algorithms. This article develops an adaptive MLMC algorithm that takes advantage of \textit{a posteriori} analysis of differential equations. In contrast to the adaptive MLMC method in \cite{KornYou_2018}, the algorithm here uses adjoint-based \textit{a posteriori} error approximations to adaptively refine meshes for new levels of the MLMC estimator and also accurately determine when the estimator has achieved a specified stopping criteria. The bias of the MLMC estimator contains a sum of the errors in samples of a quantity of interest (QoI). As such, an adjoint-based error analysis is employed to accurately estimate these errors, and thus the bias of the MLMC estimator. This analysis provides a decomposition of the bias in order to adaptively refine meshes as new levels of the estimator are created. By adaptively refining meshes, this robust MLMC algorithm yields a more cost effective estimate than MLMC with uniformly refined meshes.

This article presents the results of analysis for both ordinary differential equation (ODE) and PDE settings. We present a classical \textit{a posteriori} error analysis which deals with `standard' QoIs that can be represented as bounded functionals of the solution that has been widely explored \cite{AinOde00, BanRan03, Bar04, BecRan01, CaoPet04, CarEstTav08, ChaEstGinTav15, ChaEstTav2019, Cha18, ChaColSha17, ChaEstGinShaTav15, ChaEstTavCarSan16, ChaEstGinTav13, ChaShaWil19, ColEstTav15, ColEstTav14, DavLeV16, EriEstHanJoh95, Est95, Est09, EstHolMik02, EstLarWil00, JohChaCarEstGinLarTav15}. The error estimation utilizes generalized Green's functions solving an adjoint problem, computable residuals of the numerical solution, and variational analysis \cite{AinOde00, BecRan01, DavLeV16, Est04, EstHolMik02, EstLarWil00}. The error estimation method is often used with a finite element or variational numerical method, however many finite difference and finite volume methods are acceptable, provided they be recast as an equivalent finite element method \cite{ChaEstGinShaTav15, ChaEstGinTav15, ColEstTav15,ColEstTav14, DelHagTro81, DelDub86, EriJohLog04, EstGinTav12, Log04,ChaEstSteTav21}. Nonlinear QoIs are typically handled by linearizing around a computed solution \cite{BecRan01, CarEstTav08, CliColHou15}. We also present the results from recent work \cite{ChaEstSteTav21,CEGST_SWE} that has extended this analysis to a certain `nonstandard' QoI which cannot be represented as a linear functional of the solution or trivially linearized.

This article begins with a review of the MLMC method for a general QoI and its mean squared error in \S \ref{sec:MLMC}. We then introduce a model PDE and model ODE along with corresponding QoIs in \S \ref{sec:Apps_QoIs}. Section \ref{sec:error_analysis} presents the adjoint-based error analysis needed to accurately compute the bias term in the  MLMC method. Different grid refinement methods based on the error estimates are presented in \S \ref{sec:adapt_refine}, and an Adaptive MLMC algorithm  is developed in \S \ref{sec:algo}. Finally, numerical  examples demonstrating the efficacy of the Adaptive MLMC algorithm are given  in \S \ref{sec:examples}.

\section{MLMC for Differential Equations}\label{sec:MLMC}
This section reviews the multi-level Monte Carlo estimator for the expected value of a random variable \cite{Giles2008,Heinrich_MLMC,giles2015}. In particular, the random variable in this discussion is a QoI related to the solution of a differential equation that depends on some random parameter(s). The particular differential equations and QoIs investigated in this article are described in more detail in \S \ref{sec:Apps_QoIs}. For now, the notation is left more general to emphasize that the MLMC estimator is applicable to any QoI that is dependent on a random parameter.

Let $Q=Q(u;w)$ be a QoI related to the solution $u$ of a differential equation that depends on a random variable $w$ belonging to the probability space $\left(\Omega, \mathcal{F}, \mathbb{P}  \right)$. Also let $Q(w^{(n)})=Q(u;w^{(n)})$, for $n=1,\dots,N$, be a large number $N$ of samples of $Q$, where $w^{(n)}$ is sampled from the same probability space $\left(\Omega, \mathcal{F}, \mathbb{P}  \right)$. Given an approximation (e.g. a numerical solution) $U$ to $u$ (i.e. $U \approx u$), let $\widehat{Q}(w^{(n)})=\widehat{Q}(U;w^{(n)})$ be an approximation of $Q(w^{(n)})$. The standard Monte Carlo method uses the $N$ samples $\widehat{Q}(w^{(n)})$ to estimate the expected value, $\mathbb{E}[Q]$, of $Q$:
\begin{align}\label{eq:MC}
    \mathbb{E}[Q] \approx \frac{1}{N}\sum_{n=1}^N\widehat{Q}(w^{(n)}).
\end{align}
There are two sources of error for the estimator \eqref{eq:MC}; error from using a finite number of samples, and discretization error from the approximation $\widehat Q \approx Q$.

The multi-level Monte Carlo method is an extension of the MC method that utilizes several different approximations $\{\widehat{Q}_{\ell}(w) = \widehat{Q}(U_{\ell};w) \}$ in order to obtain a more cost-efficient approximation of the expected value.
More rigorously, 
let $\{\widehat{Q}_{\ell}(w) = \widehat{Q}(U_{\ell};w) \}_{\ell=0}^{L-1}$ be a collection of approximate QoI where the accuracy of the approximations $U_{\ell}$ increases with $\ell$. At each level $\ell$, $N_{\ell}$ samples of the random variable $w$ are taken from the probability space $\left(\Omega, \mathcal{F}, \mathbb{P}  \right)$ and denoted as $\{ w_{\ell}^{(n)} \}_{n=1}^{N_{\ell}}$. The MLMC estimator using $L$ levels is constructed by expanding the expected value of the most accurate estimator as  \cite{Giles2008}
\begin{align}
    \mathbb{E}[Q] \approx \mathbb{E}[\widehat{Q}_{L-1}] = \mathbb{E}[\widehat{Q}_{0}] + \sum_{\ell=1}^{L-1}\mathbb{E}[\widehat{Q}_{\ell} -\widehat{Q}_{\ell-1} ]. \label{eq:expval_expand}
\end{align}
Using the standard MC method to obtain the expected values on the right side of \eqref{eq:expval_expand} gives the L-level MLMC estimator for the expected value of $Q$:
\begin{align}
    \mathbb{E}[Q] \approx \widehat Q_{L,\{N_{\ell} \}}^{ML} 
    &= \frac{1}{N_0}\sum_{n=1}^{N_0} \widehat Q_0(w_0^{(n)}) + 
    \sum_{\ell=1}^{L-1}\left\{\frac{1}{N_{\ell}}\sum_{n=1}^{N_{\ell}}\left( \widehat Q_{\ell}(w_{\ell}^{(n)}) - \widehat Q_{\ell-1}(w_{\ell}^{(n)})\right)  \right\}, \nonumber \\
    &=\sum_{\ell=0}^{L-1}\left\{\frac{1}{N_{\ell}}\sum_{n=1}^{N_{\ell}}\left( \widehat Q_{\ell}(w_{\ell}^{(n)}) - \widehat Q_{\ell-1}(w_{\ell}^{(n)})\right)  \right\},\label{eq:MLMC_est}
\end{align}
where $\widehat Q_{-1}\equiv 0$. One sample on the $\ell-th$ level is
\begin{equation}\label{eq:sample}
    Y_{\ell}(w_{\ell}^n)=\left( \widehat Q_{\ell}(w_{\ell}^{(n)}) - \widehat Q_{\ell-1}(w_{\ell}^{(n)})\right),
\end{equation}
and requires two different approximations of the QoI. Both come with the same sample of the random parameter $w$, but they are obtained using the different approximations $\widehat Q_{\ell}$ and $\widehat Q_{\ell-1}$.


\subsection{Mean Squared Error for MLMC method}
The Mean Squared Error (MSE) of the MLMC estimator is \cite{giles2015}
\begin{align}\label{eq:MSE}
    MSE &=\mathbb{E}\left[\left(\widehat{Q}^{ML}_{ \{N_{\ell}\},L  } - Q \right)^2 \right]= \mathbb{V}\left[ \widehat{Q}^{ML}_{ \{N_{\ell}\},L  }  \right] + \left( \mathbb{E}\left[\widehat{Q}^{ML}_{ \{N_{\ell}\},L  }  - Q  \right] \right)^2 \nonumber  \\
    &= \sum_{\ell=0}^{L-1} \frac{1}{N_{\ell}}\mathbb{V}\left[\widehat{Q}_{\ell} - \widehat{Q}_{\ell-1} \right] + \left( \mathbb{E}\left[\widehat{Q}_{L-1} - Q  \right] \right)^2, 
 \end{align}
 where the variance $\mathbb{V}$ of a random variable $X$ is
 \begin{equation}
     \mathbb{V}[X]= \mathbb{E}\left[\left(X-\mathbb{E}[X]  \right)^2  \right].
 \end{equation}
 The first term in \eqref{eq:MSE} is the variance, which is decomposed into contributions from each level of the multi-level estimator. The second term is the squared bias, which only depends on the highest level. If a tolerance $\epsilon$ is desired for the MSE, it is sufficient that
 \begin{equation}
     \mathbb{V}\left[ \widehat{Q}^{ML}_{ \{N_{\ell}\},L  }  \right]< \epsilon/2 \quad \text{ and } \quad \left( \mathbb{E}\left[\widehat{Q}^{ML}_{ \{N_{\ell}\},L  }  - Q  \right] \right)^2 < \epsilon/2.
 \end{equation}
 
 The variance at level $\ell$ and the bias can be approximated using a finite number of samples \cite{giles2015}
\begin{align}
    \mathbb{V}\left[\widehat{Q}_{\ell} - \widehat{Q}_{\ell-1} \right] &\approx \frac{1}{N_{\ell}-1}\sum_{n=1}^{N_{\ell}} \left( \widehat{Q}_{\ell}(w_{\ell}^{(n)}) - \widehat{Q}_{\ell-1}(w_{\ell}^{(n)}) -  \frac{1}{N_{\ell}}\sum_{k=1}^{N_{\ell}} \left(\widehat{Q}_{\ell}(w_{\ell}^{(k)}) - \widehat{Q}_{\ell-1}(w_{\ell}^{(k)}) \right) \right)^2,\label{eq:var_l} \\
    \mathbb{E}\left[\widehat{Q}_{L-1} - Q  \right] &\approx \frac{1}{N_{L-1}}\sum_{n=1}^{N_{L-1}}\left(\widehat{Q}_{L-1}(w_{L-1}^{(n)}) - Q(w_{L-1}^{(n)})\right).\label{eq:bias}
\end{align}
With the approximations \eqref{eq:var_l} and \eqref{eq:bias}, the MSE is approximated as
 \begin{align}\label{eq:MSE_est}
         MSE&\approx \sum_{\ell=0}^{L-1} \frac{1}{N_{\ell}} \left[  \frac{1}{N_{\ell}-1}\sum_{n=1}^{N_{\ell}} \left( \widehat{Q}_{\ell}(w_{\ell}^{(n)}) - \widehat{Q}_{\ell-1}(w_{\ell}^{(n)}) -  \frac{1}{N_{\ell}}\sum_{k=1}^{N_{\ell}} \left(\widehat{Q}_{\ell}(w_{\ell}^{(k)}) - \widehat{Q}_{\ell-1}(w_{\ell}^{(k)}) \right) \right)^2 \right] + \nonumber \\
    & +\left(\frac{1}{N_{L-1}}\sum_{n=1}^{N_{L-1}}\left(\widehat{Q}_{L-1}(w_{L-1}^{(n)}) - Q(w_{L-1}^{(n)})\right)\right)^2. 
 \end{align}

The squared-bias term in \eqref{eq:MSE_est} is not directly computable as it requires a true sample of the QoI, $Q(w_{L-1}^{(n)})$; formulas to estimate the bias for different QoIs are presented in \S \ref{sec:error_analysis}. This estimate for the MSE shows that the bias depends on the error of the highest level of the MLMC estimator. In order to effectively lower the bias, the estimator $\widehat{Q}_{L-1}$ would need to become more accurate, i.e. a new level would have to be introduced. This new highest level generally does not require a large number of samples\cite{giles2015,Heinrich_MLMC}.

The MSE estimate \eqref{eq:MSE_est} also shows that taking more samples on a given level (increasing some $N_{\ell}$) decreases the overall variance of the MLMC estimator. However, taking more samples increases the cost of the estimator so a balance must be met in order to decrease the variance without increasing the cost significantly.
Let the variance on level $\ell$ be denoted as $V_{\ell}=\mathbb{V}\left[\widehat{Q}_{\ell} - \widehat{Q}_{\ell-1} \right]$ and let the cost of taking one sample on the $\ell$-th level be $C_{\ell}$. Using the method of Lagrange multipliers to minimize overall cost, $\sum N_{\ell}C_{\ell} $, under the constraint for the total variance $\sum \frac{V_{\ell}}{N_{\ell}} < \frac{1}{2}\epsilon$ gives the optimal number of samples to take on level $\ell$ to be \cite{giles2015}
\begin{equation}\label{eq:N_opt}
    N_{\ell,opt} =\left \lceil \frac{2}{\epsilon}\sqrt{\frac{V_{\ell}}{C_{\ell}}}\sum_{k=0}^{L-1}\sqrt{\frac{V_k}{C_k}} \right \rceil.
\end{equation}
Note that in practice only an approximation of the variances $V_{\ell}$ are available and thus only an approximate $N_{\ell,opt}$ can be computed. If this yields an under-approximation of $N_{\ell,opt}$, the MSE of the MLMC estimator will be larger than the desired tolerance.

\section{Applications and Quantities of Interest}\label{sec:Apps_QoIs}
The MLMC method as described in \S \ref{sec:MLMC} is applicable to any QoI that depends on a random parameter. This section describes different differential equations and three specific QoIs analyzed in this article and on which  numerical experiments are carried out in \S \ref{sec:examples}. First, a generic stationary PDE is presented alongside a standard QoI expressible as a linear functional. Then, generic  ODEs with initial conditions (or initial value problems) are discussed with a standard QoI and a non-standard QoI which is not expressible as a linear functional. The $L^2(\Omega)$ inner-product is denoted as the inner-product pairing $\left( \alpha,\beta \right)$, while the Euclidean inner-product is denoted by the dot-product $\alpha \cdot \beta$.

\subsection{Stationary PDEs}
Consider a partial differential equation of the form
\begin{align}\label{eq:diff_eq}
\begin{cases}
    \mathcal{D}u &= f, \hspace{.5cm} x\in \Omega,\\
    u &=g, \hspace{.5cm} x \in \partial \Omega.
\end{cases}
\end{align}
over a domain $\Omega \subset \mathbb{R}^d $ with differential operator $\mathcal{D}$  $= \mathcal{D}(w)$,
and continuous functions $f=f(u,x;w)$ and $g=g(x;w)$, any of which may depend on a random variable $w$. The solution $u=u(x;w)$ of \eqref{eq:diff_eq} will therefore depend on the random parameter $w$.
For example, in the steady-state advection-diffusion equation
\begin{equation}\label{eq:PDE_example}
    \nabla^2u + b \cdot {\nabla u} = f,
\end{equation}
the operator is $\mathcal{D}u=\nabla^2u + b \cdot {\nabla u} $, where the vector $b\in \mathbb{R}^d$ may have random components $b_i=b_i(w)$ for some $i$, $1\leq i\leq d$.

\subsubsection{QoIs for Stationary PDEs}
The form of QoI used in this article for the stationary equations is considered to be a `standard' QoI because it can be expressed as a linear functional of the solution:
\begin{align}
    Q_S(u) = \left(u,\psi\right) = \int_{\Omega}u\cdot\psi {\rm d}\Omega  ,\label{eq:S_qoi}
\end{align}
for some function $\psi\in L^2(\Omega)$ with $\psi(x)=0$ for $x\in \partial\Omega$. This QoI is used to represent a weighted average of $u$ over a smaller area inside $\Omega$. The adjoint-based error representation for this QoI is described in Theorem \ref{thm:standard_error}.

\subsection{ODEs}
This article also considers ordinary differential equations of the form
\begin{align}\label{eq:diff_eq_time}
\begin{cases}
    \frac{du}{dt} &= f, \hspace{.5cm} t\in (0,T],\\
    u(0)&=u_0,
\end{cases}
\end{align}
where $f=f(u,t;w)\in\mathbb{R}^d$ is a Lipschitz continuous function, and  $u(t)\in \mathbb{R}^d$ is the solution to the ODE. Here, the function $f$ or the initial condition $u_0= u_0(w)$ may depend on a random variable $w$. For example consider an ODE of the form (an example of this form is studied in \S \ref{sec:harm}),
\begin{equation}\label{eq:generic_lin_ode}
    \frac{d{u}}{dt}=\tilde f(t)-A{u},
\end{equation}
where $A$ is a matrix with random components $A_{i,j}=A_{i,j}(w)$ for some pair(s) $i,j$, $1\leq i,j \leq 2$.

\subsubsection{QoIs for ODEs}\label{sec:QoIs_time}
Two forms of QoI are used for ODEs. The first is a standard QoI similar to \eqref{eq:S_qoi} and involves the solution evaluated at a particular time $t^*$:
\begin{equation}\label{eq:S_qoi_time}
    Q_S(u) =u(t^*)\cdot\psi,
\end{equation}
for some $\psi\in \mathbb{R}^d$. For example, the QoI \eqref{eq:S_qoi_time} could represent the average value of the components of $u$ at the final time $t^*=T$. The error representation for this QoI is presented in Theorem \ref{thm:standard_error_time}.

\begin{figure}[ht]
\centering
\subfloat[Nonstandard QoI for the harmonic oscillator in \S\ref{sec:harm_nonstand}, yielding the fifth occurrence of $u_1=0$.]{
\includegraphics[width=5cm]{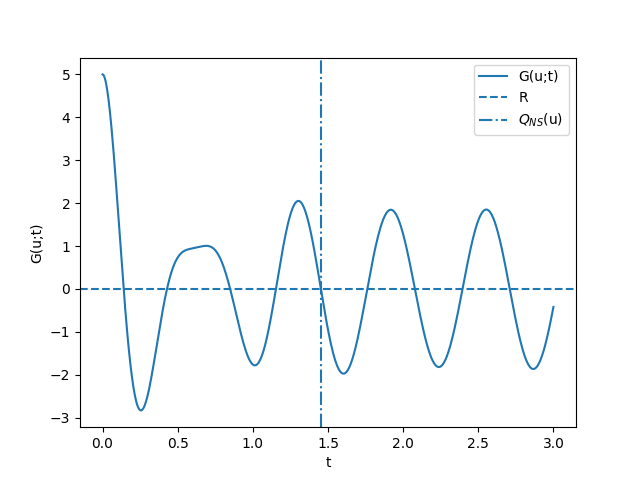}
 \label{fig:NSQoI_harm}}
 \quad
\subfloat[Accumulated error $E_k$ and corresponding meso-regions for level $\ell=0$ of example in \S \ref{sec:harm_nonstand}.]{
\includegraphics[width=5cm]{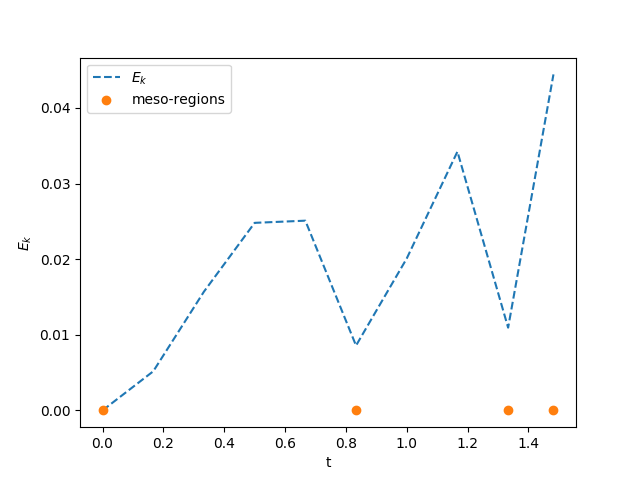}
\label{fig:accum_err}}
\caption{}
\end{figure}

The second form of QoI for time-dependent equations is the nonstandard QoI which represents the time at which a certain event related to the solution of \eqref{eq:diff_eq_time} occurs \cite{CEGST_SWE,ChaEstSteTav21}. To express this quantity, let $G(u;t)$ be a linear functional of $u$ which is implicitly dependent on $t$ through $u$, and expressed in terms of the inner product,
\begin{equation} \label{eq:G}
G(u; t) =u(t)\cdot\psi,
\end{equation}
for some $\psi\in\mathbb{R}^d$. Also, let $R$ be a chosen threshold value and assume that there are one or more times $t^\star$ during the interval $(0,T]$ for which $G(u; (t^\star)) = R$. Define the time $H(u,\hat t)$ for fixed $G$ and $R$ as
\begin{equation} \label{eq:defn_H}
H(u,\hat t) = \min_{t \in (\hat t, T]} \arg ( G(u; t) = R ).
\end{equation}
where $\hat t$ is specified in order to obtain different occurrences of the event $G(u; t)=R$. For example, if $\hat t =0$ then $H(u,\hat t) = H(u,0)$ denotes the time of the first occurrence of the event $G(u; t) = R$. If $H(u,0) \leq \hat t < H(u,H(u,0))$ then $H(u,\hat t)$ gives the time of the second occurrence of $G(u; t)=R$. The nonstandard quantity of interest $Q_{NS}(u)$ for fixed $\hat t$ is then defined as \cite{CEGST_SWE}
\begin{equation} \label{eq:QoI_time_to_event}
Q_{NS}(u) = H(u,\hat t).
\end{equation}
A depiction of this type of QoI is shown in Figure \ref{fig:NSQoI_harm}, where $\hat t$ has been chosen to obtain the fifth occurrence of $G(u; t) = R$; see \S \ref{sec:harm_nonstand} for details. The error representation for this non-standard QoI is shown in Theorem \ref{thm:NS_error_rep}.
\section{Adjoint-based Error Analysis}\label{sec:error_analysis}
The computation of the bias in MLMC, given in \eqref{eq:bias}, requires estimation of the error in the QoI for each sampled parameter $w_{L-1}^{(n)}$. This section presents {a posteriori} analysis for estimating this error for a given (or fixed) parameter $w$ and derives error formulation of the different QoIs discussed in \S \ref{sec:Apps_QoIs}. Since the sampled parameter $w$ is fixed throughout this section, we omit mentioning the dependence of solutions on the parameter for simplicity of notation.
The error representations are derived using variational analysis and rely on solutions to certain adjoint problems. First, the variational forms of the two differential equations \eqref{eq:diff_eq} and \eqref{eq:diff_eq_time} are presented. Then, the error representations and the associated adjoint problems are discussed for the three QoIs \eqref{eq:S_qoi}, \eqref{eq:S_qoi_time}, and \eqref{eq:QoI_time_to_event}.

\subsection{Variational Forms of Differential Equations}\label{sec:var_forms}
Let $W$ be a Hilbert space over $\Omega$ with inner-product $\left(\alpha,\beta \right)_W$ and the induced W-norm $\| \alpha \|_W$, for $\alpha,\beta \in W$.

\noindent \textbf{Stationary PDEs}: The variational form of \eqref{eq:diff_eq} can be expressed as: Find $u \in W$ such that 
\begin{align}\label{eq:weak_eq}
    \left(\mathcal{D}_1u,\mathcal{D}_2v\right) = \left(f,v\right), \hspace{1cm} \forall v \in W, 
\end{align}
with $L^2(\Omega)$ inner-product $\left(a,b\right) = \int_{\Omega} a\cdot b \ {\rm d}\Omega$. For example, for the PDE \eqref{eq:PDE_example}, $W=H_0^1(\Omega)$. The operators $\mathcal{D}_1,\mathcal{D}_2$ are linear differential operators in $x$, such that $\mathcal{D}_2^*\mathcal{D}_1 u = \mathcal{D}u$. The operator $\mathcal{D}_2^*$ is the formal adjoint of $\mathcal{D}_2$ which satisfies the property that
\begin{equation}\label{eq:formal_adjoint}
    \left(u,\mathcal{D}_2v\right) = \left(\mathcal{D}_2^*u,v\right), \hspace{0.5cm} \text{ for all } u \in W, v \in W.
\end{equation}
The differential operator $\mathcal{D}_1$ may be of lower order than $\mathcal{D}$. Although, if a solution $u$ to \eqref{eq:weak_eq} is sufficiently differentiable, it will also satisfy \eqref{eq:diff_eq}.

\noindent \textbf{ODEs}: The variational form of \eqref{eq:diff_eq_time} is: Find $u \in \left(H^1((0,T))\right)^d$ 
such that
\begin{equation} \label{eq:diff_eq_time_weak}
\int_0^T \frac{du}{dt}\cdot v {\rm d}t  = \int_0^T  f \cdot v {\rm d}t, \quad \forall v \in \left(L^2(0,T)\right)^d.
\end{equation}
The variational equations \eqref{eq:weak_eq} and \eqref{eq:diff_eq_time_weak} are utilized to derive the error representations presented in the next section.
\subsection{Error Formulations and Adjoint Problems}\label{sec:QoI_errs}
This section provides adjoint-based error formulations for each of the QoI \eqref{eq:S_qoi}, \eqref{eq:S_qoi_time}, and \eqref{eq:QoI_time_to_event}. The formulations rely on solutions to certain adjoint problems, which are also described in this section. 
In this section, $U=U(w^{(n)})$ refers to a numerical solution to either \eqref{eq:diff_eq} or \eqref{eq:diff_eq_time}, which will be obvious from context, for a single sample of the random parameter $w^{(n)}$.

\subsubsection{Error for Standard QoIs}
The error representations for two standard QoIs are presented. Theorem \ref{thm:standard_error} applies to QoIs of form \eqref{eq:S_qoi} related to the solution of a time-independent PDE \eqref{eq:diff_eq}. Theorem \ref{thm:standard_error_time} is applicable to QoIs of form \eqref{eq:S_qoi_time} corresponding to the solution of a time-dependent differential equation \eqref{eq:diff_eq_time}.

\noindent \textbf{Stationary PDEs:}
\begin{theorem}\label{thm:standard_error}
Given a numerical solution $U(x)$ to \eqref{eq:diff_eq}, let $e(x)=u(x)-U(x)$. The error, $\left(e,\psi\right)$, in the quantity \eqref{eq:S_qoi} is given by
\begin{equation}\label{eq:standard_error_rep}
    \left(e,\psi\right) =  \left(\mathcal{D}_1U,\mathcal{D}_2\phi\right) - \left(f,\phi\right)
\end{equation}
where $\phi(x)$ is the solution of the adjoint problem
\begin{align}
\begin{cases}
    \mathcal{D}^*(\phi(x))+\psi&=0, \hspace{.5cm} x \in \Omega, \\
    \phi(x)&=0, \hspace{.5cm} x \in \partial\Omega.    
\end{cases}
\end{align}
\end{theorem}

The linear differential operator $\mathcal{D}^*$ is the corresponding adjoint operator of $\mathcal{D}$ in \eqref{eq:diff_eq} with property \eqref{eq:formal_adjoint}. The proof of Theorem \ref{thm:standard_error} is standard; see \cite{eehj_actanum_95}.

The error representation \eqref{eq:standard_error_rep} provides a way to compute the terms in the bias \eqref{eq:bias}. Moreover, taking advantage of the linearity of the integral with respect to the domain produces a decomposition of the error into contributions from different regions of the domain. Let $\mathcal{T}$ be a partition of the domain $\Omega$. Then 
\begin{align}
    \left(e,\psi\right) &=  \left(\mathcal{D}_1U,\mathcal{D}_2\phi\right) - \left(f,\phi\right) = \int_{\Omega}\left(\left[\mathcal{D}_1U\cdot \mathcal{D}_2\phi\right] -f\cdot\phi \right){\rm d}\Omega, \nonumber \\
    &=\sum_{\tau \in \mathcal{T}}\int_{\tau}\left(\left[\mathcal{D}_1U\cdot \mathcal{D}_2\phi\right] -f\cdot\phi\right) {\rm d}\Omega = \sum_{\tau \in \mathcal{T}} e_{\tau} ,\label{eq:err_decomp}
\end{align}
where the error contribution, $e_{\tau}$, from a particular simplex $\tau$ is given by the integral
\begin{equation}\label{eq:error_contributions}
    e_{\tau}=\int_{\tau}\left(\left[\mathcal{D}_1U\cdot \mathcal{D}_2\phi\right] -f\cdot\phi\right) {\rm d}\Omega.
\end{equation}


\noindent\textbf{ODEs}:
\begin{theorem}\label{thm:standard_error_time}
Given a numerical solution $U(t)$ to \eqref{eq:diff_eq_time}, let $e(t)=u(t)-U(t)$. The error, $e(t^*)\cdot\psi$, in the quantity \eqref{eq:S_qoi_time} at a time $t^*\in(0,T]$ is given by
\begin{equation}\label{eq:standard_error_rep_time}
   Q_S(u)  - Q_S(U)  =  e(t^*)\cdot\psi =  \int_0^{t^*}\left[f\cdot\phi  -\frac{dU}{dt}\cdot\phi \right] {\rm d}t,
\end{equation}
where $\phi(t)$ is the solution of the adjoint problem
\begin{align}\label{eq:std_time_adj}
\begin{cases}
    -\frac{d\phi}{dt}&= \overline{A_{u,U}(t)}^T\phi, \hspace{.5cm}, t \in[0,t^*), \\
    \phi(t^*) &= \psi.
\end{cases}
\end{align}
where
\begin{equation}
    \overline{A_{u,U}(t)} = \int_0^1 \frac{df}{dz}(z,t) {\rm d}s
\end{equation}
with $z=su+(1-s)U$.
\end{theorem}
The proof of Theorem \ref{thm:standard_error_time} is standard; see \cite{Est04}. Note that the initial condition for the adjoint problem is at time $t=t^*$ and the equation is solved backwards in time.
Since $\overline{A_{u,U}}$ requires the true solution $u$ of \eqref{eq:diff_eq_time}, it is often approximated by $\overline{A_{u,U}} \approx \nabla_uf(U,x,t)$. With this, the right side of the adjoint equation \eqref{eq:std_time_adj} is approximated by
\begin{align}
    \overline{A_{u,U}(t)}^T\phi \approx \left( \nabla_uf(U,x,t) \right)^{\top} \phi(x,t).
\end{align}

Again, taking advantage of linearity in the integral in \eqref{eq:standard_error_rep_time} provides a decomposition of the error. If the domain $(0,t^*)$ is partitioned into $N^*$ sub-intervals with endpoints $\{t_0,t_1,\dots,t_{N^*} \}$, a decomposition of the error \eqref{eq:standard_error_rep_time} is
\begin{align}\label{eq:error_decomp_time}
     e(t^*)\cdot\psi = \sum_{i=0}^{N^*-1} \int_{t_i}^{t_{i+1}}\left[f\cdot\phi-\frac{dU}{dt}\cdot\phi \right] {\rm d}t = \sum_{i=0}^{N^*-1} e_{(t_i,t_{i+1})},
\end{align}
with error contributions, $e_{(t_i,t_{i+1})}$, given by the integrals
\begin{equation}\label{eq:e_I}
    e_{(t_i,t_{i+1})} = \int_{t_i}^{t_{i+1}}\left[f\cdot\phi  -\frac{dU}{dt}\cdot\phi \right] {\rm d}t.
\end{equation}

\subsubsection{Nonstandard QoI}\label{sec:err_NSQoI}

For a chosen threshold $R$ and fixed value of $\hat t$ in \eqref{eq:defn_H}, denote the true value of the QoI as
\begin{equation}
\label{def:t_t}
  t_t:=Q(u).
\end{equation}
Let $U(x,t)$ be an approximation of the solution satisfying \eqref{eq:diff_eq_time_weak} and let
\begin{equation}
\label{def:t_c}
t_c := Q(U)
\end{equation}
be the corresponding computed (or approximated) QoI. Denote the error in the approximate solution as $e(x,t) = u(x,t) - U(x,t)$. The error representation for the nonstandard QoI \eqref{eq:QoI_time_to_event} is given in Theorem \ref{thm:NS_error_rep}.

\begin{theorem}\label{thm:NS_error_rep}
Let $u$ be the true solution to \eqref{eq:diff_eq_time_weak} and $t_t$ be the true NSQoI \eqref{eq:QoI_time_to_event}. Also let $U$ and $t_c$ be a numerically computed solution and its corresponding computed NSQoI, respectively. For a chosen functional $G$ and threshold value $R$, the error in the NSQoI \eqref{eq:QoI_time_to_event} is given as
\begin{align}\label{eq:err_rep_NS}
    e_Q &= t_t-t_c  = \frac{e(t_c)\cdot\psi  + \mathcal{R}_1 }{f(t_c)\cdot\psi  +  e(t_c)\cdot \nabla_u[\psi \cdot f] +\mathcal{R}_2}.
\end{align}
Where $e(t) = u(t) - U(t)$ and the remainders satisfy $\mathcal{R}_1 = \mathcal{O}\left((t_t-t_c)^2\right)$ and $\mathcal{R}_2 = \mathcal{O}\left(|u-U|^2\right)$.
\end{theorem}
The proof of Theorem \ref{thm:NS_error_rep} relies on Taylor's Theorem to expand around known quantities; see \cite{ChaEstSteTav21}.

In the error representation \eqref{eq:err_rep_NS}, the magnitude of the remainders decrease faster than other terms as the numerical solution becomes more accurate. Setting  $ \mathcal{R}_1  \approx 0$ and $\mathcal{R}_2 \approx 0$  gives the approximation

\begin{equation}\label{eq:approx1}
    t_t-t_c  \approx  
    \frac{e(t_c)\cdot\psi  }{f(t_c)\cdot\psi  +  e(t_c)\cdot \nabla_u[\psi \cdot f]  }. 
\end{equation}
The terms containing the unknown error $e(t_c)$ are estimated using Theorem \ref{thm:standard_error_time}, requiring one adjoint problem per term.

The numerator of \eqref{eq:approx1} is decomposed in a similar fashion as \eqref{eq:error_decomp_time}, with $t^*=t_c$. This gives the decomposition of the error in the non-standard QoI as
\begin{align}\label{eq:NS_err_decomp}
    t_t-t_c \approx \frac{1}{\mathcal{C}}\sum_{i=0}^{N^*-1}e_{(t_i,t_{i+1})}
\end{align}
where $\mathcal{C}= f(t_c)\cdot\psi  +  e(t_c)\cdot \nabla_u[\psi \cdot f] $ is the non-decomposed denominator of the error estimate.

The three error decompositions, \eqref{eq:err_decomp}, \eqref{eq:error_decomp_time}, and \eqref{eq:NS_err_decomp}, are used to adaptively refine meshes when creating new levels of the MLMC estimator, as is discussed in the next section.

\section{Apaptive Multi-level Monte Carlo Algorithm}\label{sec:MLMC_Algo}

This section develops an MLMC algorithm that uses the error decompositions from \S \ref{sec:error_analysis} to adaptively create meshes as new levels are added to the estimator. First, two grid refinement methods are presented. 


\subsection{Creating New Levels for MLMC}\label{sec:adapt_refine}
 There are many refinement methods to choose from when creating the grid for a new level. MLMC algorithms usually employ standard uniform refinement \cite{giles2015}. We take advantage of the form of the error decompositions \eqref{eq:err_decomp}, \eqref{eq:error_decomp_time}, and \eqref{eq:NS_err_decomp} in order to adaptively create new grids and in turn create a more efficient MLMC algorithm. The rest of this section describes the two adaptive refinement methods used in this article. Since these refinement methods are applied to individual samples, we also discuss how each refinement method is adapted to deal with multiple samples.

\subsubsection{Dual Weighted Residual Refinement}
For the dual-weighted-residual (DWR) refinement method, the regions, $\tau$ (used in the error decompositions; see \S \ref{sec:QoI_errs}), corresponding to the largest contributions of error are refined by a given factor \cite{BecRan01,BanRan03}. This type of grid refinement relies on the decomposition of the error \eqref{eq:err_decomp} and some criteria to determine whether a region should be refined or not. Possible criteria include refining any region that corresponds to an error larger than some tolerance, or refining a certain number of regions that contribute the most error.

When using this refinement method in an MLMC algorithm, we want to combine data from multiple samples to determine which regions to refine. We do this by determining which regions should be refined for each individual sample, and then refining \textit{all} of these regions. 

More precisely, let $\mathcal{T}_h=\{\tau_1,\tau_2,...,\tau_M\}$ be a simplicial decomposition of a domain $\Omega$, where $h$ denotes the maximum diameter of the elements of $\mathcal{T}_h$. Let $U(w^{(1)})=U(x;w^{(1)})$ and $U(w^{(2)})=U(x;w^{(2)})$ be two numerical solutions corresponding to different samples of the random parameter $w$. Let $\mathcal{J}=\{\tau_{J_1},\tau_{J_2},...,\tau_{J_{\widehat M}}\}$ be the set of $\widehat M < M$ regions to be refined based on the error decomposition corresponding to the numerical solution $U(w^{(1)})$. Also let $\mathcal{K}=\{\tau_{K_1},\tau_{K_2},...,\tau_{K_{\widetilde M}}\}$ be the similar set corresponding to $U(w^{(2)})$. We then refine all of the regions in the union of both sets, $\mathcal{J} \cup \mathcal{K}$.



\subsubsection{Meso-scale Refinement}
In time-dependent problems, a meso-scale refinement of the time domain may be used in order to preserve cancellations of the error over large sections of the domain~\cite{ChaEstTavCarSan16}. This refinement method considers the accumulation of the error over sub-intervals and forms ``meso-scale regions'', i.e. regions of maximal cancellation in the error. The meso-scale regions are each uniformly refined, using different scaling factors on different regions. Since uniform refinement usually preserves cancellations, the locations of these minima are preserved but the error decreases. Examples presented in this article follow the ``Allocation of fixed resources'' algorithm from \cite{ChaEstTavCarSan16}, which is described below.

Given a numerical solution $U$ of an ODE \eqref{eq:diff_eq}, computed over a temporal grid $\mathcal{T} =\{I_1,I_2,\dots,I_{\widetilde N} \}$, define the accumulated contributions to the error as
\begin{equation}
   E_k =\left|\sum_{i=1}^k e_{I_i} \right|,
\end{equation}
where the $e_{I_i}$ are define in \eqref{eq:e_I}. Meso-scale regions, or meso-regions, are determined by finding the global minimum of the accumulated error after the initial increase of accumulated error. This is then repeated, starting at the minimum, until the end of the temporal domain is reached. See Figure \ref{fig:accum_err} for an example of an $E_k$ and the corresponding meso-regions.

Let $P$ be the number of meso-scale regions and  ${\widetilde N}_i$ be the number of time-steps in the $i$-th meso-scale region.  If the $i$-th meso-scale region starts at interval $I_{p}$ and ends at interval $I_r$, then $\widetilde N_i=r-p+1$. Also let $\mathcal{E}_i = (E_r-E_p)$ denote the error accumulated over the $i$-th meso-scale region. Assume that this error accumulated over the $i$-th meso-scale region satisfies
\begin{equation}\label{eq:meso_errs}
    \mathcal{E}_i= \frac{c_i}{{\widetilde N}_i^{q}},
\end{equation}
where $c_i$ is some positive constant and $q$ is determined by the order of the numerical method. The goal of the meso-scale refinement is to create a new grid with $\widehat N$ sub-intervals that minimizes the total error. If the total number of intervals $\widehat N$ is fixed, the total error is minimized if
\begin{equation}\label{eq:N_criteria}
    \frac{c_i}{{\widehat N}_i^{q+1}}= K, \ \ \ \ \forall i,
\end{equation}
for some constant $K$. The $\widehat N_is$ are obtained as follows:
First obtain all $c_i$ from \eqref{eq:meso_errs}. Then, rearranging \eqref{eq:N_criteria} to get an expression for $\widehat N_i$ and using the fact that
\begin{equation}
    \sum_{i=1}^P \widehat N_i = \widehat N,
\end{equation}
compute the value of $K$ as
\begin{equation}
    K=\left[\frac{1}{\widehat N}\sum_{i=1}^P\left( c_i^{1/(q+1)}\right) \right]^{q+1}.
\end{equation}
Finally, obtain the $\widehat N_i$s from \eqref{eq:N_criteria}.

It is not clear how to modify the meso-scale refinement to combine information from multiple samples. Instead of combining the samples, we use the refinement that corresponds to the single sample that had the largest computed absolute error in the QoI,
with the added condition that regions are never unrefined.

\begin{figure}[ht]
\centering
\subfloat[Meso-regions]{
\includegraphics[width=5cm]{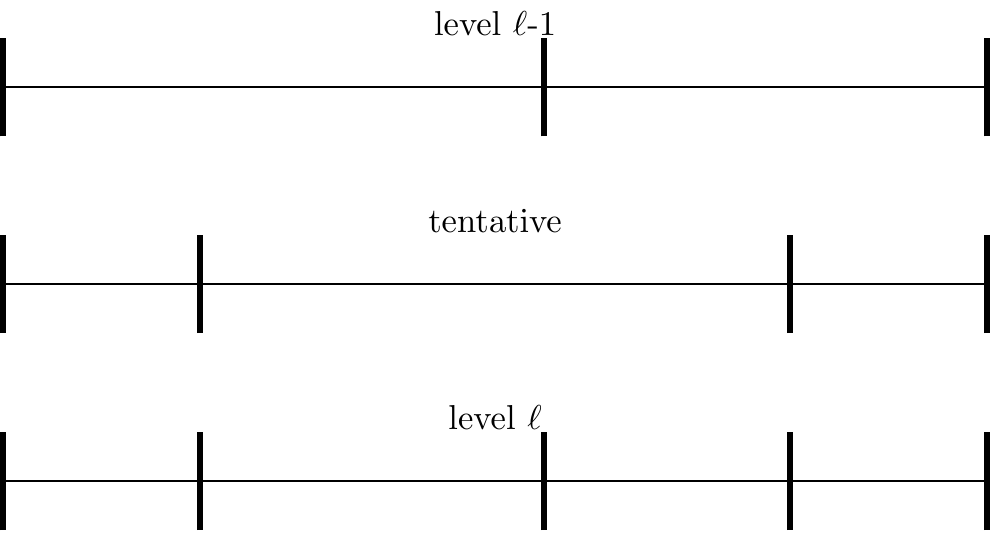}
\label{fig:meso_regions}}
\hspace{1.0in}
\subfloat[Meso-regions and sub-intervals]{
\includegraphics[width=5cm]{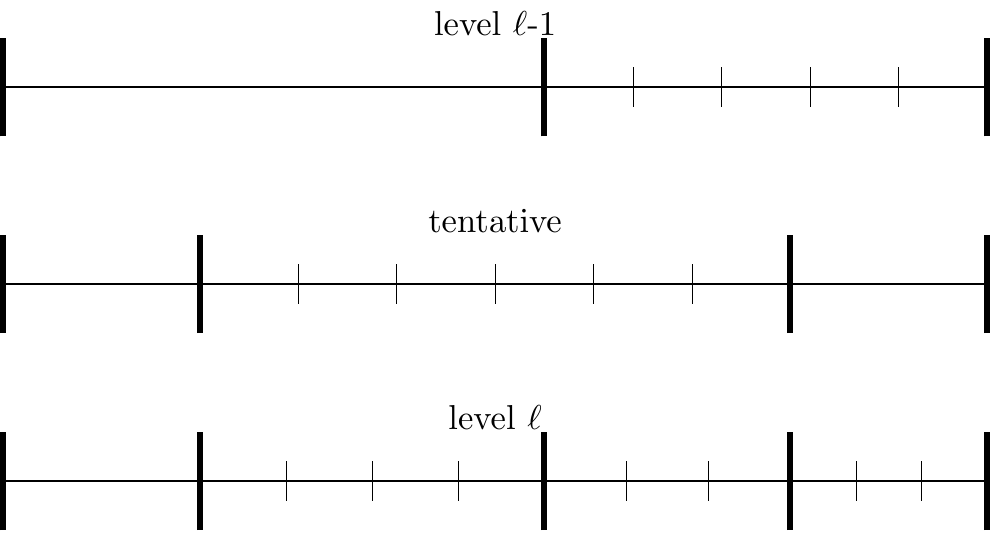}
\label{fig:meso_example}}
\caption{From top to bottom: level $\ell-1$, tentative, and level $\ell$ grid. The meso-regions are marked by thicker and longer lines, while the intervals within each meso-region are marked by lighter and shorter lines.}
\end{figure}

This is accomplished as follows. First, a tentative level $\ell$ grid is created with the meso-regions and the $\widehat N_i$s as described above. Then, a common refinement of the meso-regions from this tentative and the level $\ell -1$ is taken, see Figure \ref{fig:meso_regions} for an illustration. For each of the regions in this common refinement, determine whichever grid ($\ell-1$ or tentative) has more nodes in that region. Finally, using that number of nodes, make a uniform sub-grid over that region. This process is explained in Figure \ref{fig:meso_example}. Here the grid at level $\ell - 1$ has two meso-scale regions, with the first meso-region having a single interval while the second meso-region has five intervals. The tentative grid (shown in the middle) has three meso-regions having one, six and one intervals respectively. The refinement, which forms the level $\ell$ grid is shown at the bottom. This grid has four meso-regions having one, four, three and three intervals respectively. We observe that all regions of level $\ell$ grid have a finer discretization than the corresponding regions in the level $\ell-1$ grid.

\subsection{Description of Adaptive Algorithm}\label{sec:algo}
The adaptive MLMC routine is given in Algorithm \ref{alg:driver}. It requires a user supplied initial mesh, an initial number of samples $N$, and a tolerance $\epsilon$ for which the goal is to achieve $MSE<\epsilon$. First $N$ samples \eqref{eq:sample} involving the QoI \eqref{eq:S_qoi} are obtained. Each of these $N$ samples requires sampling the random parameter $w$, solving \eqref{eq:diff_eq} numerically over the initial provided mesh, and using that solution to obtain the QoI. Along with each sample, the contributions \eqref{eq:error_contributions} to the error in the computed QoI are obtained and the sample variance \eqref{eq:var_l} and sample bias \eqref{eq:bias} are computed. With the sample variance, the optimal number of samples \eqref{eq:N_opt} is computed. If the optimal number of samples is larger than $N$, more samples \eqref{eq:sample} involving the QoI \eqref{eq:S_qoi} are taken and the bias and variance are updated to include the new data.

Next, the bias is compared to the desired tolerance. If the tolerance has not been achieved, a new level is added. The mesh for the new level is constructed using an adaptive refinement method as discussed in \S \ref{sec:adapt_refine}, based off the error contributions \eqref{eq:error_contributions}. $N$ samples \eqref{eq:sample} involving the QoI \eqref{eq:S_qoi}
are obtained on this new level, along with the error contributions \eqref{eq:error_contributions},variance \eqref{eq:var_l}, and bias \eqref{eq:bias}. The optimal number of samples \eqref{eq:N_opt} is computed for \textit{all} levels. Extra samples are taken on each level as needed and the variances are updated. The bias is also updated to include all samples taken on the highest level. This process is repeated, adding as many levels as needed, until the bias has reached the desired tolerance. The MSE is then computed via \eqref{eq:MSE_est}.


\begin{algorithm}
\caption{Adaptive MLMC driver routine}\label{alg:driver}
\begin{algorithmic}
    \Require{$N$,init$\_$mesh,$\epsilon$}
    \State Set L=1 and $\ell=0$
    \State Compute $N$ samples \eqref{eq:sample} involving the QoI \eqref{eq:S_qoi}, error contributions \eqref{eq:error_contributions}, variance \eqref{eq:var_l}, and bias \eqref{eq:bias}.\;
    \State Find optimal number of samples, $N_{0,opt}$, using \eqref{eq:N_opt}.\;
    \State Compute $N_{0,opt}-N$ new samples \eqref{eq:sample} involving the QoI \eqref{eq:S_qoi}, update variance \eqref{eq:var_l} and bias \eqref{eq:bias} to include all $N_{0,opt}$ samples.\;
    \While{$bias^2>\epsilon/2$}
    \State Add new level, L=L+1 and $\ell=\ell+1$\;
    \State Create new mesh.\;
    \State Compute $N$ samples \eqref{eq:sample} involving the QoI \eqref{eq:S_qoi}, error contributions \eqref{eq:error_contributions}, variance \eqref{eq:var_l}, and bias \eqref{eq:bias}.\;
    \For{$\widehat\ell$ in 0,...,L-1}
    \State Find optimal number of samples \eqref{eq:N_opt} for level $\widehat\ell$. \;
    \State Take extra samples \eqref{eq:sample} as needed. Update variance \eqref{eq:var_l} to include new samples.\;
    \EndFor
    \State Compute bias \eqref{eq:bias}.\;
    \EndWhile
    \State Compute MSE using \eqref{eq:MSE_est}\;
    \State Compute MLMC estimator, $\widehat Q_{L,\{N_{\ell} \}}^{ML}$, using \eqref{eq:MLMC_est} and the MSE using \eqref{eq:MSE_est}\;
    \State \textbf{Output:} $\widehat Q_{L,\{N_{\ell} \}}^{ML}$, MSE
\end{algorithmic}
\end{algorithm}

\begin{remark}
When creating new levels, the initial number of samples taken does not have to be the same number $N$ as for the lowest level. If the same large $N$ is used for every level, starting new levels will be much more expensive than starting previous levels and it becomes more likely that unnecessary samples are taken on high levels. In our examples below, we keep the cost of starting a level relatively fixed for the first three levels, after which we keep the number of initial samples constant. This maintains a reasonable computational cost while also taking enough samples for the variance estimates to be meaningful.
\end{remark}

\begin{remark}
Algorithm \ref{alg:driver} details how to obtain an estimated expected value in a standard QoI \eqref{eq:S_qoi} corresponding to the solution of the differential equation \eqref{eq:diff_eq}, although the algorithm is applicable to any differential equation and QoI, provided an error decomposition is available.
\end{remark}

\section{Numerical Results}\label{sec:examples}

This section begins by discussing the finite element method (FEM) used in the numerical examples. For each problem, we present the adjoint equation needed for the error analysis of \S \ref{sec:error_analysis}. Examples include standard QoIs for a harmonic oscillator and a stationary advection-diffusion equation, as well as non-standard QoIs for a harmonic oscillator, the Lorenz equations, and the two-body problem. The MLMC algorithm is applied to each problem using different refinement methods with results provided for comparison.

\subsection{Numerical Methods: Continuous Galerkin FEM}
 The error analysis presented in \S \ref{sec:error_analysis} is well-suited for numerical solutions obtained from variational methods. The analysis also applies to other numerical methods, provided they be cast to an equivalent variational form. The examples in this section use the first order, continuous Galerkin finite element method.

Continuous Galerkin finite element method are constructed using the standard continuous Lagrangian finite element spaces. Let $\mathcal{T}_h$ be a simplicial decomposition of $\Omega$, where h is the
maximum diameter of the elements of $\mathcal{T}_h$, such that the union of the elements of $\mathcal{T}_h$ is $\Omega$, and the intersection of any two elements is either a common edge, a node, or is empty. The degree $q$ continuous Lagrange finite element space is then defined as
\begin{equation}
    \mathcal{P}_h^q = \left\{v\in C(\Omega) : \forall \tau \in \mathcal{T}_h, \left.v\right|_{\tau} \in \mathcal{P}^q(\tau)  \right\},
\end{equation}
where $\mathcal{P}^q(\tau) $ is the space of polynomials of degree at most $q$ defined on the element $\tau$.

\noindent \textbf{Stationary PDEs}: The degree $q$ continuous Galerkin finite element method, with respect to $\mathcal{T}_h$, for the PDE \eqref{eq:diff_eq} is: Find $U\in \mathcal{P}_h^q$ such that
\begin{align}
    \left(\mathcal{D}_1U,\mathcal{D}_2v\right) = \left(f,v\right), \hspace{1cm} \forall v \in \mathcal{P}^q_h.
\end{align}

\noindent \textbf{ODEs}: For the differential equation \eqref{eq:diff_eq_time}, the domain is the interval $(0,T]$ and $\mathcal{T}_h$ is a set of sub-intervals $\{(t_{n},t_{n+1}) \}_{n=0}^{N^*-1}$ with $t_{n+1}-t_n \leq h$. The continuous Galerkin method is: Find $U\in \mathcal{P}_h^q$ such that the restriction of $U$ to any $(t_{n},t_{n+1}) \in \mathcal{T}_h$ satisfies
\begin{align}
    \int_{t_n}^{t_{n+1}} \frac{dU}{dt}\cdot v {\rm d}t = \int_{t_n}^{t_{n+1}} f\cdot v {\rm d}t , \hspace{1cm} \forall v \in \mathcal{P}^{q-1}(t_n,t_{n+1}).
\end{align}

For example, for the ODE \eqref{eq:generic_lin_ode} with a domain decomposed into sub-intervals $\mathcal{T}_h$, the continuous Galerkin method with $q=1$ is: Find $U\in \mathcal{P}_h^1$ such that the restriction of $U$ to any $(t_{n},t_{n+1}) \in \mathcal{T}_h$ satisfies
\begin{align}
    \int_{t_n}^{t_{n+1}} \frac{dU}{dt}\cdot v {\rm d}t = \int_{t_n}^{t_{n+1}} (\tilde f - AU)\cdot v {\rm d}t , \hspace{1cm} \forall v \in \mathbb{R}^2.
\end{align} 
Here, the solution $U$ is continuous and piece-wise linear with respect to the partition $\mathcal{T}_h$.

\subsection{Harmonic Oscillator}\label{sec:harm}
Consider the harmonic oscillator
\begin{equation*}
\ddot{\omega}=
  -\frac{k}{m}\omega-\frac{c}{m}\dot{\omega}+\frac{F_0}{m}\cos(\nu t +\theta_d), \ t \in (0,3], \quad \omega(0)=5, \ \dot{\omega}(0)=0.
\end{equation*}
with deterministic parameters $c=1, \ F_0=50, \ \theta_d=0, \ \nu=10,$
and where $k$ and $m$ are random parameters.
Rewriting as a system of first-order ODEs, $\dot{u}+Au=\tilde f(t)$, gives
\begin{equation}\label{eq:harm_osc}
\begin{pmatrix}
\dot{u_1}(t) \\  \dot{u_2}(t)
\end{pmatrix}
+
\begin{pmatrix}
0    & -1 \\
k/m  & 1/m
\end{pmatrix}
\begin{pmatrix}
u_1(t) \\  u_2(t)
\end{pmatrix}
=
\begin{pmatrix}
0 \\  50/m * \cos(10 t)
\end{pmatrix},
\ t \in (0,3].
\end{equation}
With initial conditions $(u_1(0),u_2(0))=(5,0)$.
We look at examples of both a standard QoI and the aforementioned non-standard QoI.

\subsubsection{Standard QoI}\label{sec:harm_stand}

Consider the oscillator \eqref{eq:harm_osc} where $k\sim N(50,2)$, a normal distribution with mean 50 and standard deviation 2, and $m\sim Unif[0.225,0.275]$, a uniform distribution between 0.225 and 0.275. The standard QoI is $Q_S(u)= [(1,0)^{\top}\cdot u(T)] = u_1(T)$, the position of the oscillator at the final time $T=3$. For this QoI, the error representation \eqref{eq:standard_error_rep_time} is
\begin{equation}
     e(T)\cdot(1,0)^{\top}  = \int_0^T \phi\cdot\left(\tilde f-\dot{U}-AU \right){\rm d}t,
\end{equation}
where $\phi$ is the solution to the adjoint problem
\begin{align}
\begin{cases}
    -\dot{\phi} +A^{\top}\phi &= 0 \hspace{.5cm} t \in [0,T),\\
    \phi(T)&=1.
\end{cases}
\end{align}

The grid for the lowest level has 27 elements, with the number of elements roughly doubling for each further level. The algorithm is run using three different grid creation methods; uniform grid creation, DWR refinement, and meso-scale refinement. In each, we start by taking 100 samples on the lowest level. When creating further levels, the second level starts with 50 samples and all further levels start with 20 samples. The tolerance for the MSE is set to $\epsilon=0.001$.

Results when using uniform grids are shown in Tables \ref{tab:harm_stand_unif} and \ref{tab:harm_stand_unif1}. Table \ref{tab:harm_stand_unif} provides details for each level of the estimator, including the number of elements used to create the grid, the relative cost per sample, the number of samples, and variance contribution \eqref{eq:var_l}. Table \ref{tab:harm_stand_unif1} gives results for the MLMC estimator, including the variance, squared bias, MSE \eqref{eq:MSE_est}, the estimated expected value \eqref{eq:MLMC_est}, and the total relative cost. The cost is computed by summing across all levels, the number of samples at each level times the cost per sample  at that level. The cost of the sample at level $\ell=0$ (i.e. the cost of computing a single value of $\widehat{Q}_0(w_{\ell}^{(n)})$, for any $\ell$) is normalized to 1.
Samples for levels $\ell>0$ require computing two values, $\widehat{Q}_{\ell}(w_{\ell}^{(n)})$ and $\widehat{Q}_{\ell-1}(w_{\ell}^{(n)})$. Hence, the cost of taking a sample on level $\ell$ is the sum of the costs of computing $\widehat{Q}_{\ell}(w_{\ell}^{(n)})$ and $\widehat{Q}_{\ell-1}(w_{\ell}^{(n)})$ relative to the cost of computing $\widehat{Q}_0(w_{\ell}^{(n)})$. For example, if $\ell=1$ and computing a single value of $\widehat{Q}_{1}(w_{1}^{(n)})$ is twice the computational cost of computing a single $\widehat{Q}_{0}(w_{\ell}^{(n)})$, then the relative cost of a sample on the $\ell=1$ level is $Cost(\widehat{Q}_{1}(w_{1}^{(n)})) + Cost(\widehat{Q}_{0}(w_{\ell}^{(n)})) = 2+1 =3$.

\begin{figure}[ht]
\centering
\subfloat[Uniform refinement grids]{
\includegraphics[width=5cm]{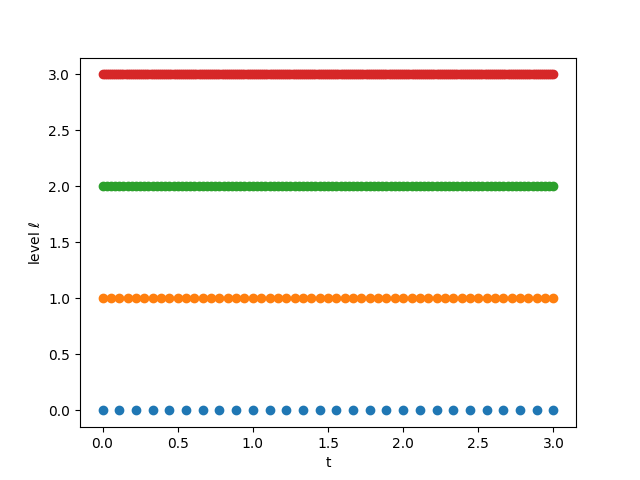}
\label{fig:osc_stand_unif}}
\subfloat[DWR refinement grids]{
\includegraphics[width=5cm]{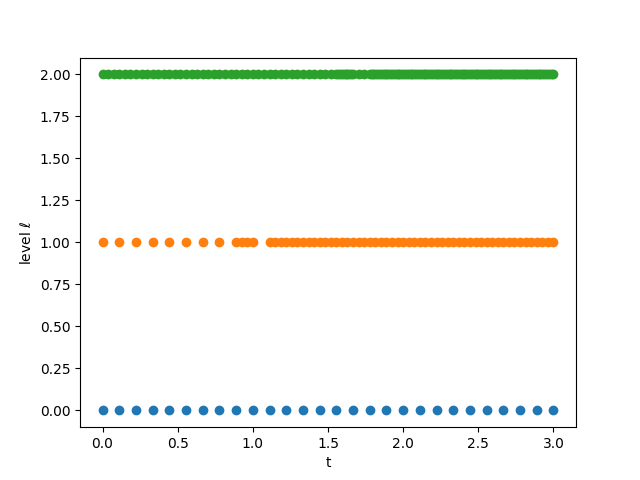}
\label{fig:osc_stand_adapt}}
\subfloat[Meso-scale refinement grids]{
\includegraphics[width=5cm]{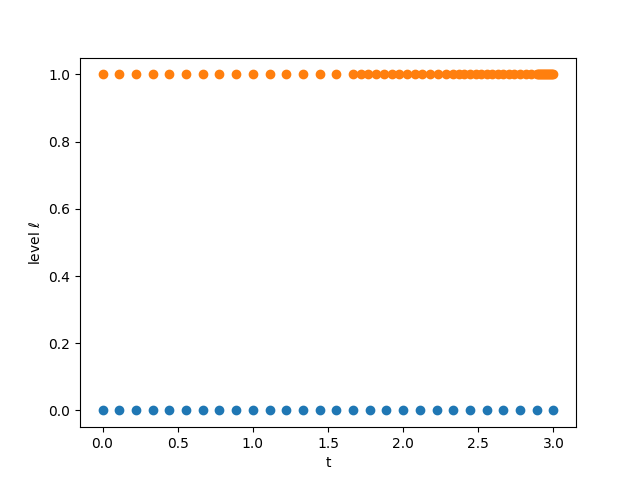}
\label{fig:osc_stand_meso}}
\caption{Grids from the different refinement methods for the example in \S \ref{sec:harm_stand}.}\label{fig:harm1}
\end{figure}

Tables \ref{tab:harm_stand_adapt} and \ref{tab:harm_stand_adapt1} show similar details when the grids are created using DWR refinement. Results for the MLMC estimator using meso-scale refinement are provided in Tables \ref{tab:harm_stand_meso} and \ref{tab:harm_stand_meso1}.
When using uniform grids, the MLMC estimator requires four levels in order to achieve tolerance in bias. The estimator using DWR refinement requires three levels. With meso-scale refinement, only two levels are required to have a bias less than tolerance. More samples are required on level $\ell=1$ when using DWR or meso-scale refinement but the overall costs are still lower compared to using uniform refinement.  Grids for each refinement method are shown in Figure \ref{fig:harm1}.
\begin{table}[h!]
    \centering
    \begin{tabular}{|c|c|c|c|c|}
         \hline
         Level & \# Elems & Cost Per Sample & \# Samples & Var. Per Level \\
         \hline
         0 & 27   & 1  & 227  & 4.32968E-4 \\
         \hline
         1 & 54  & 3  & 52 & 7.36927E-4 \\
         \hline
         2 & 108  & 6   &  20 & 1.90041E-6 \\
         \hline
         3 & 216  & 12   &  20 & 3.28090E-7 \\
         \hline
    \end{tabular}
    \caption{Results for each level of estimator in example from \S \ref{sec:harm_stand} with $\epsilon=.001$. New grids are obtained via uniform refinement.}
    \label{tab:harm_stand_unif}
\end{table}
\begin{table}[h!]
    \centering
    \begin{tabular}{|c|c|c|c|c|}
    \hline
        Tot. Var. & Squared Bias & MSE & Est. Exp. Val & Tot. Cost \\
    \hline
       0.00117 & 2.82413E-5  & 0.00120  & -0.38276 & 743 \\
    \hline
    \end{tabular}
    \caption{Results of MLMC estimator in example from \S \ref{sec:harm_stand} with $\epsilon=.001$. New grids are obtained via uniform refinement.}
    \label{tab:harm_stand_unif1}
\end{table}
\begin{table}[h!]
    \centering
    \begin{tabular}{|c|c|c|c|c|}
         \hline
         Level & \# Elems & Cost Per Sample & \# Samples & Var. Per Level \\
         \hline
         0 & 27   & 1  & 227  & 4.32968E-4 \\
         \hline
         1 &  63  & 3.333  & 75  & 4387715E-4  \\
         \hline
         2 &  153 & 8   &  20 & 1.36574E-6  \\
         \hline
    \end{tabular}
    \caption{Results for each level of estimator in example from \S \ref{sec:harm_stand} with $\epsilon=.001$. New grids are obtained via DWR refinement where the 50\% largest contributions to the error are refined by a factor of 3.}
    \label{tab:harm_stand_adapt}
\end{table}
\begin{table}[h!]
    \centering
    \begin{tabular}{|c|c|c|c|c|}
    \hline
        Tot. Var. & Squared Bias & MSE & Est. Exp. Val & Tot. Cost \\
    \hline
        9.22049E-4 & 8.02556E-5  &  0.001002 & -0.40173 & 637  \\
    \hline
    \end{tabular}
    \caption{Results of MLMC estimator in example from \S \ref{sec:harm_stand} with $\epsilon=.001$. New grids are obtained via DWR refinement where the 50\% largest contributions to the error are refined by a factor of 3.}
    \label{tab:harm_stand_adapt1}
\end{table}
\begin{table}[h!]
    \centering
    \begin{tabular}{|c|c|c|c|c|}
         \hline
         Level & \# Elems & Cost Per Sample & \# Samples & Var Per Level \\
         \hline
         0 & 27  & 1  & 227  & 4.32968E-4 \\
         \hline
         1 & 59  & 3.185  & 65  & 8.27785E-4   \\
         \hline
    \end{tabular}
    \caption{Results for each level of estimator in example from \S \ref{sec:harm_stand} with $\epsilon=.001$. New grids are obtained via meso-scale refinement.}
    \label{tab:harm_stand_meso}
\end{table}
\begin{table}[h!]
    \centering
    \begin{tabular}{|c|c|c|c|c|}
    \hline
        Tot. Var & Squared Bias & MSE & Est. Exp. Val &Tot. Cost\\
    \hline
     0.00126 & 1.80659E-5 & 0.00127 & -0.3816 & 434.03 \\
    \hline
    \end{tabular}
    \caption{Results of MLMC estimator in example from \S \ref{sec:harm_stand} with $\epsilon=.001$. New grids are obtained via meso-scale refinement.}
    \label{tab:harm_stand_meso1}
\end{table}

\subsubsection{Nonstandard QoI}\label{sec:harm_nonstand}
Using the same  equation and deterministic parameters as in \S \ref{sec:harm_stand}, now let $k\sim N(50,1)$ and  $m\sim Unif(0.235,0.265)$. Let the nonstandard QoI to be the time of the 5th occurrence of $u_1=0$. More precisely, set $\psi=(1,0)^\top$ in \eqref{eq:G}. Also, for the function $H$ in \eqref{eq:defn_H}, set $R=0$ and choose a $\hat t$ between the fourth and fifth occurrence of $u_1=0$; see Figure \ref{fig:NSQoI_harm}.
The error representation \eqref{eq:err_rep_NS} becomes
\begin{equation}
        t_t-t_c  \approx
    \frac{\psi\cdot e( t_c) )}{ A^{\top}\psi\cdot U(t_c)- \psi\cdot\tilde f(t_c)  + A^{\top}\psi\cdot  e(t_c)  }. 
\end{equation}
Using Theorem \ref{thm:standard_error_time}, the two error terms are given as
\begin{align}
    \psi\cdot e(t_c)  &= \int_0^{t_c}\left[\tilde f\cdot\phi_1  -\frac{dU}{dt}\cdot\phi_1  - AU\cdot\phi_1  \right] {\rm d}t\\
    A^{\top}\psi\cdot e(t_c) &= \int_0^{t_c}\left[\tilde f\cdot\phi_2 -\frac{dU}{dt}\cdot\phi_2 - AU\cdot\phi_2  \right] {\rm d}t,
\end{align}
where $\phi_1$ and $\phi_2$ are the solutions to the adjoint problems:
\begin{align}
\begin{cases}
    -\dot{\phi_1} +A^{\top}\phi_1 &= 0 \hspace{.5cm} t \in [0,t_c),\\
    \phi_1(t_c)&=\psi,
\end{cases}, \qquad
\begin{cases}
    -\dot{\phi_2} +A^{\top}\phi_2 &= 0 \hspace{.5cm} t \in [0,t_c),\\
    \phi_2(t_c)&=A^\top\psi,
\end{cases}.
\end{align}
All examples in this section start with 100 samples of a numerical solution obtained over a uniform grid with 18 sub-intervals. The second level starts with 50 samples and all further levels with 20 samples. The grids used for levels beyond the first are obtained from the different creation methods as discussed in \S \ref{sec:adapt_refine}. 

\begin{figure}[h!]
\centering
\subfloat[Uniform refinement grids]{
\includegraphics[width=5cm]{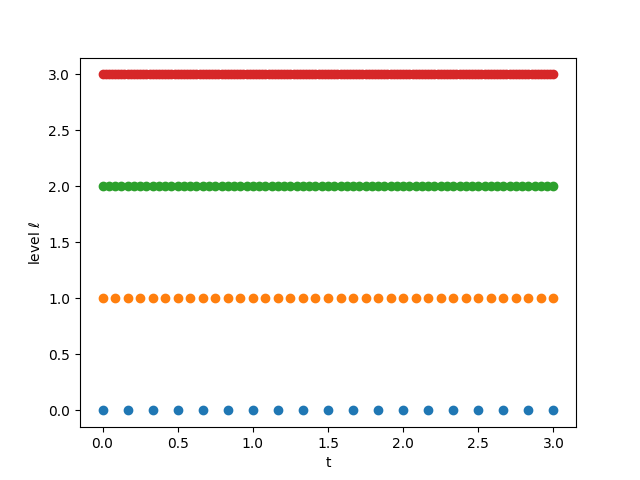}
\label{fig:osc_ns_unif}}
\subfloat[DWR refinement grids]{
\includegraphics[width=5cm]{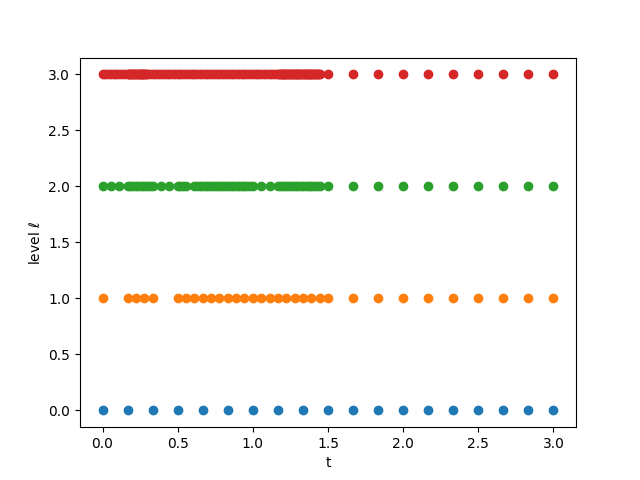}
\label{fig:osc_ns_adapt}}
\subfloat[Meso-scale refinement grids]{
\includegraphics[width=5cm]{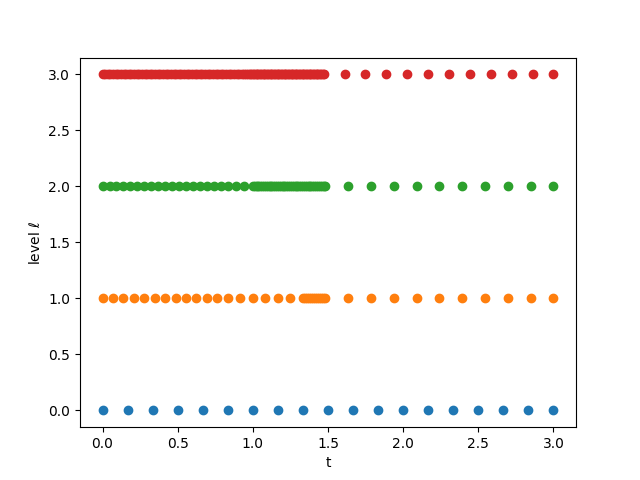}
\label{fig:osc_ns_meso}}
\caption{Grids from the different refinement methods for the example in \S \ref{sec:harm_nonstand}.}
\label{fig:osc_nonstand_grids}
\end{figure}

Tables \ref{tab:harm_nonstand_unif} and \ref{tab:harm_nonstand_unif1} provide results for the levels and overall MLMC estimator when using uniform grids. Results for the MLMC estimator using DWR refinement are shown in Tables \ref{tab:harm_nonstand_adapt} and \ref{tab:harm_nonstand_adapt1}. Tables \ref{tab:harm_nonstand_meso} and \ref{tab:harm_nonstand_meso1} give results for the estimator when grids are created using meso-scale refinement. The grids for different levels of each method can be seen in Figure \ref{fig:osc_nonstand_grids}. In all three cases, the initial number of samples is enough to achieve tolerance for the variance of the estimators. The different grid creation methods perform similarly, all requiring four levels and the same number of samples. The two adaptive methods are slightly more cost efficient because the highest level grids use less nodes than in the uniform grid method.
\begin{table}[h!]
    \centering
    \begin{tabular}{|c|c|c|c|c|}
         \hline
         Level  & \# Elems & Cost Per Sample & \# Samples & Var. Per Level \\
         \hline
         0      &  18     & 1       & 100  & 1.04977E-6  \\
         \hline
         1      &  36     &  3   & 50  & 1.44317E-6 \\
         \hline
         2      &  72     &  6   & 20  & 5.22219E-6 \\
         \hline
         3      &  144     &  12   & 20  & 2.11021E-9 \\
         \hline
    \end{tabular}
    \caption{Results for each level of estimator in example from \S \ref{sec:harm_nonstand} with $\epsilon=1E-5$. New grids are obtained via uniform refinement.}
    \label{tab:harm_nonstand_unif}
\end{table}
\begin{table}[h!]
    \centering
    \begin{tabular}{|c|c|c|c|c|}
    \hline
        Tot. Var. & Squared Bias & MSE & Est. Exp. Val & Tot. Cost \\
    \hline
      3.01727E-6 & 2.02532E-6  & 5.04259E-6 & 1.45701 & 610 \\
    \hline
    \end{tabular}
    \caption{Results of MLMC estimator in example from \S \ref{sec:harm_nonstand} with $\epsilon=1E-5$. New grids are obtained via uniform refinement.}
    \label{tab:harm_nonstand_unif1}
\end{table}
\begin{table}[h!]
    \centering
    \begin{tabular}{|c|c|c|c|c|}
         \hline
         Level  & \# Elems & Cost Per Sample & \# Samples & Var. Per Level \\
         \hline
         0      &  18     & 1       & 100  &  1.04977E-6 \\
         \hline
         1      &  34     &  2.888   & 50  & 2.17999E-6 \\
         \hline
         2      &   72    &  5.888   & 20  & 1.54181E-6 \\
         \hline
         3      &   140    &  11.777   & 20  & 2.25459E-9 \\
         \hline
    \end{tabular}
    \caption{Results for each level of estimator in example from \S \ref{sec:harm_nonstand} with $\epsilon=1E-5$.  New grids are obtained via DWR refinement where the 50\% largest contributions to the error are refined by a factor of 3.}
    \label{tab:harm_nonstand_adapt}
\end{table}
\begin{table}[h!]
    \centering
    \begin{tabular}{|c|c|c|c|c|}
    \hline
        Tot. Var. & Squared Bias & MSE & Est. Exp. Val & Tot. Cost\\
    \hline
      4.77384E-6 & 1.73603E-7 & 4.94744E-6  & 1.45620 & 597.77 \\
    \hline
    \end{tabular}
    \caption{Results of MLMC estimator in example from \S \ref{sec:harm_nonstand} with $\epsilon=1E-5$. New grids are obtained via DWR refinement where the 50\% largest contributions to the error are refined by a factor of 3.}
    \label{tab:harm_nonstand_adapt1}
\end{table}
\begin{table}[h!]
    \centering
    \begin{tabular}{|c|c|c|c|c|}
         \hline
         Level  & \# Elems & Cost Per Sample & \# Samples & Var. Per Level \\
         \hline
         0      &  18     & 1       & 100  & 1.0497E-6 \\
         \hline
         1      &  39     & 3.166    & 50  & 1.06958E-6 \\
         \hline
         2      &  70     &   6.055  & 20  & 1.23273E-7 \\
         \hline
          3      &  121     &  10.611   & 20  & 3.83435E-9 \\
         \hline
    \end{tabular}
    \caption{Results for each level of estimator in example from \S \ref{sec:harm_nonstand} with $\epsilon=1E-5$. New grids are obtained via meso-scale refinement.}
    \label{tab:harm_nonstand_meso}
\end{table}
\begin{table}[h!]
    \centering
    \begin{tabular}{|c|c|c|c|c|}
    \hline
        Tot. Var. & Squared Bias & MSE & Est. Exp. Val & Tot. Cost \\
    \hline
      2.24646E-6 & 2.11483E-7 & 2.45794E-6 & 1.45597 & 591.66 \\
    \hline
    \end{tabular}
    \caption{Results of MLMC estimator in example from \S \ref{sec:harm_nonstand} with $\epsilon=1E-5$. New grids are obtained via meso-scale refinement.}
    \label{tab:harm_nonstand_meso1}
\end{table}
\subsection{Lorenz Equations}\label{sec:lor_NS}
Consider the (nonlinear) Lorenz system,
\begin{equation}\label{eq:Lorenz}
\left.
\begin{gathered}\begin{aligned}
    \dot{u}_1 &= \sigma(u_2 - u_1), \\
    \dot{u}_2 &= r u_1 - u_2 - u_1 u_3, \\
    \dot{u}_3 &= u_1 u_2 - b u_3,
\end{aligned}\end{gathered}
\; \right \} \;
t \in (0, 2] \qquad \hbox{with} \qquad
\left \{
\begin{gathered}\begin{aligned}
    u_1(0)&=\theta, \\
    u_2(0) &= 0, \\
    u_3(0) &= 24,
\end{aligned}\end{gathered}
\; \right. \;
\end{equation}
and set $\sigma =10, r=28, \text{ and } b= \frac{8}{3}$. The initial condition $\theta$ is a random variable $\theta \sim Unif(0,2]$.

The nonstandard QoI is the time of the 2nd occurrence of $u_1 = 3$. That is, set $\psi=(1,0,0)^\top$, in \eqref{eq:G}. In \eqref{eq:defn_H} set $R=3$ and choose a $\hat t$ between the first and second occurrence of $u_1=3$.

\begin{figure}[h!]
\centering
\subfloat[Uniform refinement grids]{
\includegraphics[width=5cm]{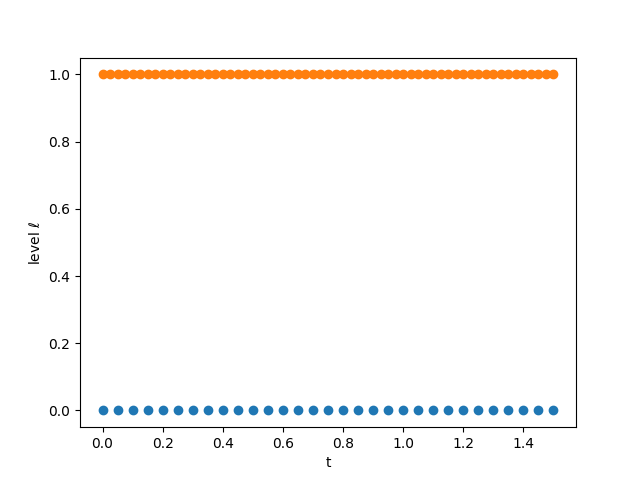}
\label{fig:lor_unif}}
\subfloat[DWR refinement grids]{
\includegraphics[width=5cm]{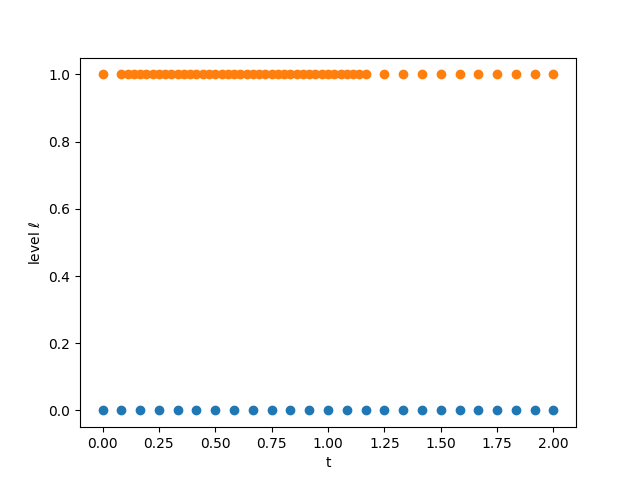}
\label{fig:lor_adapt}}
\subfloat[Meso-scale refinement grids]{
\includegraphics[width=5cm]{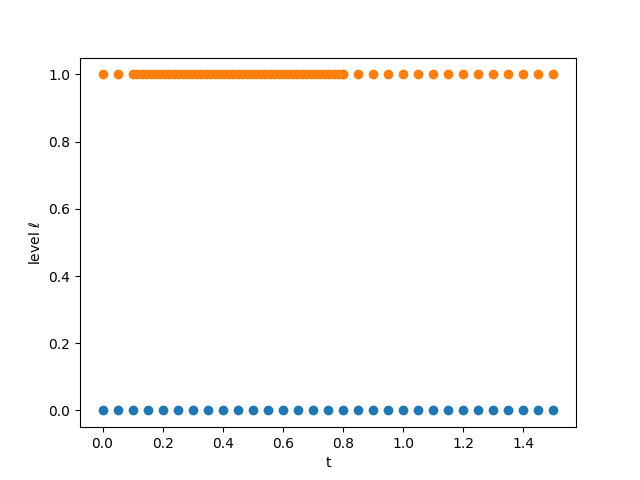}
\label{fig:lor_meso}}
\caption{Grids from the different refinement methods for the example in \S \ref{sec:lor_NS}}
\label{fig:lor_nonstand_grids}
\end{figure}

The error representation \eqref{eq:err_rep_NS} becomes
\begin{equation}\label{eq:nsqoi_ode}
        t_t-t_c  \approx
    \frac{\psi\cdot e( t_c) }{  \left(\nabla_uf(t_c)\right)^\top\psi\cdot U(t_c)- \psi\cdot f(t_c)  +  \left(\nabla_uf(t_c)\right)^\top\psi\cdot e(t_c)  }. 
\end{equation}
Using Theorem \ref{thm:standard_error_time}, the two error terms are given as
\begin{align}\label{eq:nsqoi_adjoints_ode}
    \psi\cdot e(t_c) &= \int_0^{t_c}\left[f\cdot\phi_1 -\frac{dU}{dt}\cdot\phi_1 \right] {\rm d}t\\
    \left(\nabla_uf(t_c)\right)^\top\psi\cdot e(t_c)) &= \int_0^{t_c}\left[f\cdot\phi_2 -\frac{dU}{dt}\cdot\phi_2  \right] {\rm d}t,
\end{align}
where $\phi_1$ and $\phi_2$ are the solutions to the adjoint problems
\begin{align}
\begin{cases}
    -\dot{\phi_1} &= \left(\nabla_uf\right)^\top\phi \hspace{.5cm} t \in [0,t_c),\\
    \phi_1(t_c)&=\psi,
\end{cases}, \qquad
\begin{cases}
    -\dot{\phi_2} &=\left( \nabla_uf\right)^\top\phi \hspace{.5cm} t \in [0,t_c),\\
    \phi_2(t_c)&=\left(\nabla_uf\right)^\top\psi,
\end{cases}.
\end{align}
The algorithm is started with 100 samples of a numerical solution obtained over a uniform grid with 15 sub-intervals. The second level starts with 50 samples and all further levels with 20 samples. The grids used for levels beyond the first are obtained from the different creation methods as discussed in \S \ref{sec:adapt_refine}. Tables \ref{tab:lor_nonstand_unif}, \ref{tab:lor_nonstand_adapt}, and \ref{tab:lor_nonstand_meso} show results for each level of the estimators when using uniform grids, DWR refinement and meso-scale refinement, respectively. Results for the MLMC estimators, using uniform grids, DWR refinement, and meso-scale refinement are shown in Tables \ref{tab:lor_nonstand_unif1}, \ref{tab:lor_nonstand_adapt1}, and \ref{tab:lor_nonstand_meso1}, respectively. The grids used in each estimator are shown in Figure \ref{fig:lor_nonstand_grids}. For this nonlinear problem, all three grid creation methods perform similarly. The different methods require two levels and use the same number of samples and hence the comparable approximate costs.
\begin{table}[h!]
    \centering
    \begin{tabular}{|c|c|c|c|c|}
         \hline
         Level & \# Elems & Cost Per Sample & \# Samples & Var. Per Level \\
         \hline
         0      &   24    & 1       & 101  & 4.59302E-5  \\
         \hline
         1      &  48     &  3   & 50  & 5.66121E-7 \\
         \hline
    \end{tabular}
    \caption{Results for each level of the estimator in example from \S \ref{sec:lor_NS} with $\epsilon=1E-4$. New grids are obtained via uniform refinement.}
    \label{tab:lor_nonstand_unif}
\end{table}
\begin{table}[h!]
    \centering
    \begin{tabular}{|c|c|c|c|c|}
    \hline
        Tot. Var. & Squared Bias & MSE & Est. Exp. Val & Tot. Cost \\
    \hline
      4.64963E-5 & 3.68484E-5  & 8.33448E-5  & 0.84945 & 251  \\
    \hline
    \end{tabular}
    \caption{Results of the MLMC estimator in example from \S \ref{sec:lor_NS} with $\epsilon=1E-4$. New grids are obtained via uniform refinement.}
    \label{tab:lor_nonstand_unif1}
\end{table}
\begin{table}[h!]
    \centering
    \begin{tabular}{|c|c|c|c|c|}
         \hline
         Level & \# Elems & Cost Per Sample & \# Samples & Var. Per Level \\
         \hline
         0      &   24    & 1       & 100  & 4.59302E-5  \\
         \hline
         1      &   50    &  3.083   & 50  & 8.13257E-7 \\
         \hline
    \end{tabular}
    \caption{Results for each level of the estimator in example from \S \ref{sec:lor_NS} with $\epsilon=1E-4$. New grids are obtained via DWR refinement where the 50\% largest contributions to the error are refined by a factor of 3.}
    \label{tab:lor_nonstand_adapt}
\end{table}
\begin{table}[h!]
    \centering
    \begin{tabular}{|c|c|c|c|c|}
    \hline
        Tot. Var. & Squared Bias & MSE & Est. Exp. Val & Tot. Cost\\
    \hline
      4.67435E-5 & 6.52347E-6 & 5.32670E-5  & 0.84501 & 254.16 \\
    \hline
    \end{tabular}
    \caption{Results of the MLMC estimator in example from \S \ref{sec:lor_NS} with $\epsilon=1E-4$. New grids are obtained via DWR refinement where the 50\% largest contributions to the error are refined by a factor of 3.}
    \label{tab:lor_nonstand_adapt1}
\end{table}
\begin{table}[h!]
    \centering
    \begin{tabular}{|c|c|c|c|c|}
         \hline
         Level & \# Elems & Cost Per Sample & \# Samples & Var. Per Level \\
         \hline
         0      &  24     & 1       & 100  & 4.59302E-5 \\
         \hline
         1      &   53    &  3.208   & 50  & 4.71161E-6 \\
         \hline
    \end{tabular}
    \caption{Results for each level of the estimator in example from \S \ref{sec:lor_NS} with $\epsilon=1E-4$. New grids are obtained via meso-scale refinement.}
    \label{tab:lor_nonstand_meso}
\end{table}
\begin{table}[h!]
    \centering
    \begin{tabular}{|c|c|c|c|c|}
    \hline
        Tot. Var. & Squared Bias & MSE & Est. Exp. Val & Tot. Cost \\
    \hline
       5.06418E-5 & 3.03186E-5 & 8.09605E-5 & 0.84931 & 260.4  \\
    \hline
    \end{tabular}
    \caption{Results of the MLMC estimator in example from \S \ref{sec:lor_NS} with $\epsilon=1E-4$. New grids are obtained via meso-scale refinement.}
    \label{tab:lor_nonstand_meso1}
\end{table}
\subsection{Two-Body Problem}\label{sec:twobody_NS}
Consider the two body problem
\begin{equation}\label{twobody_NS}
\left.
\begin{gathered}\begin{aligned}
\dot{u}_1 &= u_3, \\
\dot{u}_2 &= u_4, \\
\dot{u}_3 &= \frac{-u_1}{(u_1^2+u_2^2)^{3/2}}, \\
\dot{u}_4 &= \frac{-u_2}{(u_1^2+u_2^2)^{3/2}},
\end{aligned}\end{gathered}
\; \right\}
\quad t \in (0,10], \quad u(0)=(0.4, 0, 0, \theta)^\top,
\end{equation}
which models a small body orbiting a much larger body in two dimensions. Here $u_1, u_2$ are the spatial coordinates of the orbiting body relative to the larger body, and $u_3, u_4$ are the respective velocities. The last component of the initial condition, $\theta$, is a random variable $\theta \sim Unif[1.97,2]$.

For the nonstandard QoI, set $\psi=\left(1,0,0,0 \right)^\top$ in \eqref{eq:G}. Also set $R=0$ in \eqref{eq:defn_H} and choose a $\hat t$ to obtain the 3rd occurrence of $u\cdot\psi=u_1=0$. The error representation \eqref{eq:err_rep_NS} for this problem is the same as \eqref{eq:nsqoi_ode}.
From Theorem \ref{thm:standard_error_time}, the two error terms are again given by \eqref{eq:nsqoi_adjoints_ode}
where $\phi_1$ and $\phi_2$ are the solutions to the adjoint problems
\begin{align}
\begin{cases}
    -\dot{\phi_1} &= \left(\nabla_uf\right)^\top\phi \hspace{.5cm} t \in [0,t_c),\\
    \phi_1(t_c)&=\psi,
\end{cases}, \qquad
\begin{cases}
    -\dot{\phi_2} &=\left( \nabla_uf\right)^\top\phi \hspace{.5cm} t \in [0,t_c),\\
    \phi_2(t_c)&=\left( \nabla_uf\right)^\top\psi,
\end{cases}.
\end{align}

The MLMC algorithm begins with 100 samples of a numerical solution obtained over a uniform grid with 40 sub-intervals. The second level starts with 50 samples and all further levels with 20 samples. The grids used for levels beyond the first are obtained from the different creation methods as discussed in \S \ref{sec:adapt_refine}.

\begin{figure}[h!]
\centering
\subfloat[Uniform refinement grids]{
\includegraphics[width=5cm]{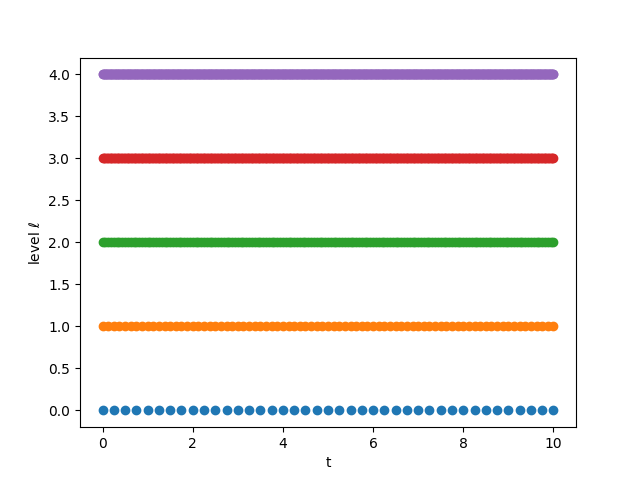}
\label{fig:tbdy_unif}}
\subfloat[DWR refinement grids]{
\includegraphics[width=5cm]{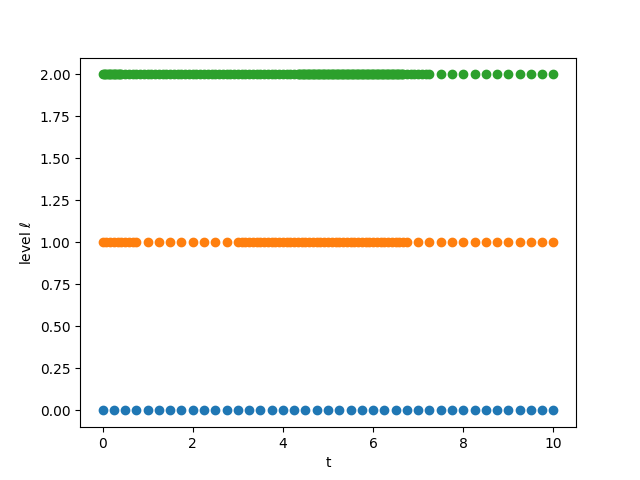}
\label{fig:tbdy_adapt}}
\subfloat[Meso-scale refinement grids]{
\includegraphics[width=5cm]{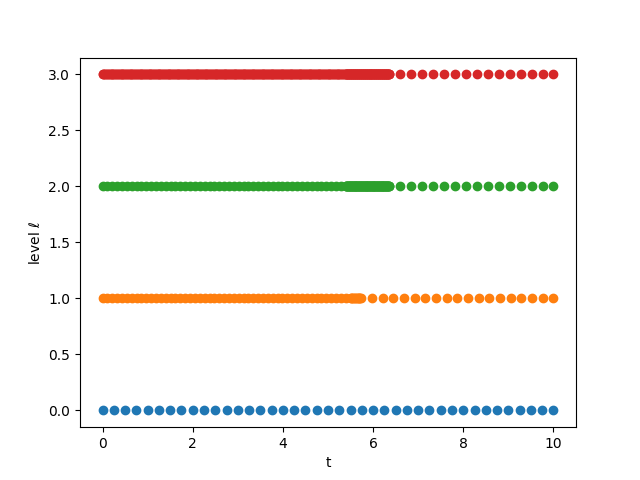}
\label{fig:tbdy_meso}}
\caption{Grids from the different refinement methods for the example in \S \ref{sec:twobody_NS}}
\label{fig:tbdy_nonstand_grids}
\end{figure}

Results for the MLMC estimator using uniform grids are shown in Tables \ref{tab:twobody_nonstand_unif} and \ref{tab:twobody_nonstand_unif1}. Tables \ref{tab:twobody_nonstand_adapt} and \ref{tab:twobody_nonstand_adapt1} give results when using DWR refinement. Results when using meso-scale refinement are provided in Tables \ref{tab:twobody_nonstand_meso} and \ref{tab:twobody_nonstand_meso1}. The grids used in each estimator are provided in Figure \ref{fig:tbdy_nonstand_grids}.

Here, the MLMC estimator using DWR refinement is the most cost efficient, requiring three levels to achieve the desired tolerance in bias. Using meso-scale refinement, the estimator requires four levels, and the estimator using uniform grids requires five levels.
Notice that the two estimators using adaptive grid creation methods do not achieve the desired tolerance for variance, due to under-sampling on some level(s). The estimator using DWR refinement is close to meeting tolerance and would only require some more samples on the lowest, cheapest level. The estimator using meso-scale refinement is further from tolerance and would require more samples on all levels, including the higher, more expensive levels. The large variance and slow convergence of the levels when using meso-scale refinement arises because our grid creation method only uses information from a single sample of the QoI.
\begin{table}[h!]
    \centering
    \begin{tabular}{|c|c|c|c|c|}
         \hline
         Level & \# Elems & Cost Per Sample & \# Samples & Var. Per Level \\
         \hline
         0      &  40     & 1   & 389   & 5.1388E-4  \\
         \hline
         1      &  80     &  3   & 50  & 1.86982E-6 \\
         \hline
         2      &  160     &  6   &  20 &  1.35918E-6 \\
         \hline
         3      &  320     &  12   &  20 & 1.20867E-7  \\
         \hline
         4      &   640    &  24   &  20 & 1.49859E-8 \\
         \hline
    \end{tabular}
    \caption{Results for each level of the estimator in example from \S \ref{sec:twobody_NS} with $\epsilon=1E-3$. New grids are obtained via uniform refinement.}
    \label{tab:twobody_nonstand_unif}
\end{table}
\begin{table}[h!]
    \centering
    \begin{tabular}{|c|c|c|c|c|}
    \hline
        Tot. Var. & Squared Bias & MSE & Est. Exp. Val & Tot. Cost \\
    \hline
      5.17251E-4 & 3.38445E-5  & 5.55109E-4  & 7.24045 & 1379  \\
    \hline
    \end{tabular}
    \caption{Results of the MLMC estimator in example from \S \ref{sec:twobody_NS} with $\epsilon=1E-3$. New grids are obtained via uniform refinement.}
    \label{tab:twobody_nonstand_unif1}
\end{table}
\begin{table}[h!]
    \centering
    \begin{tabular}{|c|c|c|c|c|}
         \hline
         Level & \# Elems & Cost Per Sample & \# Samples & Var. Per Level \\
         \hline
         0      &   40    & 1   & 255  &  6.99600E-4 \\
         \hline
         1      &   74    &  2.85   & 50  & 4.15496E-6 \\
         \hline
         2      &    160  &  5.85   & 20  & 8.88944E-7 \\
         \hline
    \end{tabular}
    \caption{Results for each level of the estimator in example from \S \ref{sec:twobody_NS} with $\epsilon=1E-3$. New grids are obtained via DWR refinement where the 50\% largest contributions to the error are refined by a factor of 3.}
    \label{tab:twobody_nonstand_adapt}
\end{table}
\begin{table}[h!]
    \centering
    \begin{tabular}{|c|c|c|c|c|}
    \hline
        Tot. Var. & Squared Bias & MSE & Est. Exp. Val & Tot. Cost\\
    \hline
      7.04644E-4 & 3.72135E-4 & 0.00107  & 7.06503 & 514.5 \\
    \hline
    \end{tabular}
    \caption{Results of the MLMC estimator in example from \S \ref{sec:twobody_NS} with $\epsilon=1E-3$. New grids are obtained via DWR refinement where the 50\% largest contributions to the error are refined by a factor of 3.}
    \label{tab:twobody_nonstand_adapt1}
\end{table}
\begin{table}[h!]
    \centering
    \begin{tabular}{|c|c|c|c|c|}
         \hline
         Level & \# Elems & Cost Per Sample & \# Samples & Var. Per Level \\
         \hline
         0      &  40   & 1       & 204  & 4.46312E-4 \\
         \hline
         1      &   81    &  3.025   & 64  & 8.85888E-4 \\
         \hline
          2      &   176    &   6.425  &  55 & 1.49708E-3 \\
         \hline
          3      &   394    &  14.25   &  20 & 1.10254E-3 \\
         \hline
    \end{tabular}
    \caption{Results for each level of the estimator in example from \S \ref{sec:twobody_NS} with $\epsilon=1E-3$. New grids are obtained via meso-scale refinement.}
    \label{tab:twobody_nonstand_meso}
\end{table}
\begin{table}[h!]
    \centering
    \begin{tabular}{|c|c|c|c|c|}
    \hline
        Tot. Var. & Squared Bias & MSE & Est. Exp. Val & Tot. Cost \\
    \hline
       2.93953E-3 & 2.18287E-4  & 0.00315 & 6.88846 & 1035.975 \\
    \hline
    \end{tabular}
    \caption{Results of the MLMC estimator in example from \S \ref{sec:twobody_NS} with $\epsilon=1E-3$. New grids are obtained via meso-scale refinement.}
    \label{tab:twobody_nonstand_meso1}
\end{table}
\subsection{Stationary Advection-Diffusion Equation}\label{sec:adv_diff}

\begin{figure}[h!]
    \centering
    \includegraphics[width=8cm]{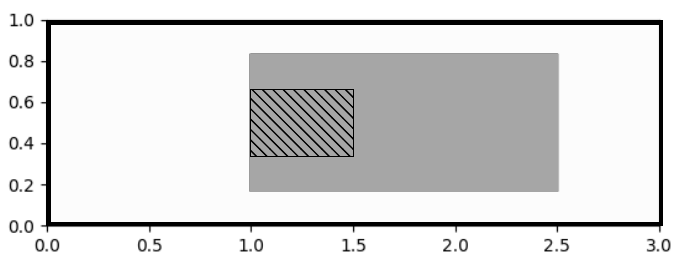}
    \caption{Illustration of overlapping supports of functions from example in Section \ref{sec:adv_diff}. The domain is $\Omega=[0,1]\times[0,3]$ outlined in black, $supp(f)=[1,2.5]\times[1/6,5/6]$ in grey, and $supp(\psi)=[1,1.5]\times[1/3,2/3]$ marked with diagonal lines.}
    \label{fig:regions_longf}
\end{figure}

\begin{figure}[h!]
\centering
\subfloat[Initial mesh, $\ell=0$]{
\includegraphics[width=7cm]{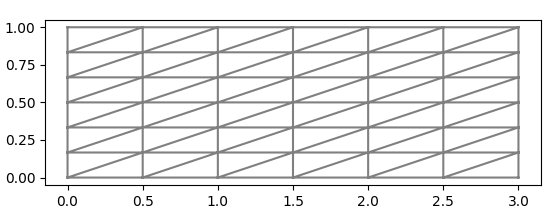}
}\\
\subfloat[Level $\ell=1$]{
\includegraphics[width=7cm]{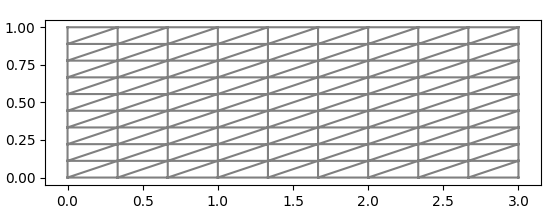}
}
\subfloat[Level $\ell=2$]{
\includegraphics[width=7cm]{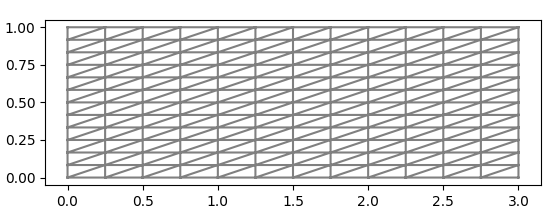}
}\\
\subfloat[Level $\ell=3$]{
\includegraphics[width=7cm]{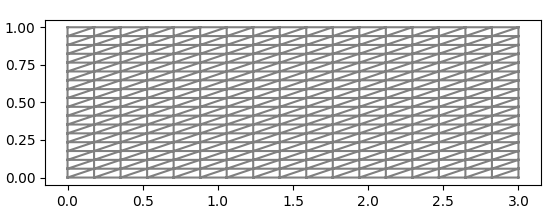}
}
\subfloat[Level $\ell=4$]{
\includegraphics[width=7cm]{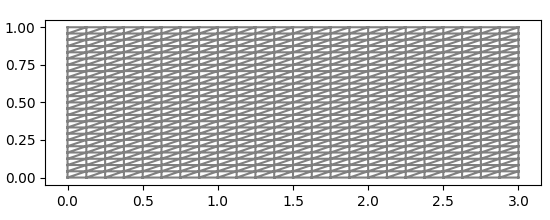}
}
\caption{Uniform meshes for example in \S \ref{sec:adv_diff}}
\label{fig:Ad_Diff_unif}
\end{figure}

\begin{figure}[h!]
\centering
\subfloat[Level $\ell=1$]{
\includegraphics[width=7cm]{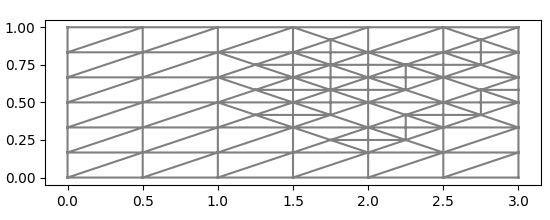}
}
\subfloat[Level $\ell=2$]{
\includegraphics[width=7cm]{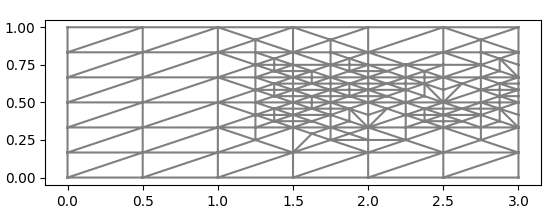}
}
\caption{Adaptively refined meshes for example in \S \ref{sec:adv_diff}}
\label{fig:Ad_Diff_adapt}
\end{figure}

Consider the equation
\begin{align}\label{eq:adv_diff}
\begin{cases}
    \nabla^2u(x) + b\cdot\nabla u(x) = f(x), \hspace{1cm} &x \in (0,3)\times(0,1), \\
    u(x)=0,\hspace{1cm} &x \in \partial\Omega.    
\end{cases}
\end{align}
The vector $b=(w,0)^{\top}$ has, as its first component, the random parameter $w=unif(1200,1600)$. The source $f$ is non-zero only over an interior region of the domain, and is given as
\begin{equation}
    f=\begin{cases}
    10000(x-1)(x-2.5)(y-1/6)(y-5/6) \ \ \ \ 1\leq x\leq 2.5, \ 1/6\leq y\leq5/6, \\
    0 \hspace{8cm} \text{else}.
    \end{cases}
\end{equation}
The weak form of \eqref{eq:adv_diff} is: Find $u \in H^1(\Omega)$ such that
\begin{equation}\label{eq:weak_p}
    -\Big(\nabla u, \nabla v \Big) + \Big(b\cdot\nabla u, v\Big) = \Big(f,v\Big), \hspace{.5cm} \forall v \in H_0^1(\Omega).
\end{equation}
The QoI is the integral of the solution over the rectangle $(1,1.5)\times(1/3,2/3)$:
\begin{equation}
    Q_S(u) = \int_{\Omega}\psi\cdot u {\rm d}\Omega, \ \ \ \text{ where } \ \ \ \psi = \begin{cases}
    1, \ \ \ \ \ \ (x,y)\in (1,1.5)\times(1/3,2/3)\\
    0, \ \ \ \ \ \ else.
    \end{cases}
\end{equation}
See Figure \ref{fig:regions_longf} for a visualization of the supports of $f$ and $\psi$. For the standard QoI \eqref{eq:S_qoi} the error representation \eqref{eq:standard_error_rep} becomes
\begin{equation}\label{eq:adv_error_rep}
    \Big(e,\psi\Big) =  - \Big(\nabla U,\nabla\phi\Big) + \Big(\nabla\cdot(bU),\phi\Big)-\Big(f,\phi\Big),
\end{equation}
where $\phi$ is the solution to the adjoint problem
\begin{align}\label{eq:adv_adj}
\begin{cases}
    \nabla^2\phi(x) - b\cdot\nabla\phi(x) + \psi(x) = 0, \hspace{1cm} & x\in \Omega,\\
    \phi(x)=0,\hspace{1cm} &x \in \partial\Omega.
\end{cases}
\end{align}
The MLMC algorithm is applied to this problem using two different mesh creation methods; uniform and DWR refinement. In both cases we aim to roughly double the number of elements when creating new levels. This is done so that each level of the different methods can be more easily compared. When using uniform meshes, we multiply the number of steps in each coordinate by $\sqrt{2}$ and round up. The uniform meshes are shown in Figure \ref{fig:Ad_Diff_unif} and results in Table \ref{tab:adv_unif} and \ref{tab:adv_unif1}. For the DWR refinement, the regions corresponding to the 25\% largest contributions to the error are refined (using DOLFIN Python's in-built refine function). The adaptively refined meshes are shown in Figure \ref{fig:Ad_Diff_adapt} and results in Tables \ref{tab:adv_adapt} and \ref{tab:adv_adapt1}. For this higher dimensional, linear problem, the MLMC estimator using DWR refinement is much more cost efficient than using uniform meshes. With uniform meshes, the estimator requires five levels to achieve the desired tolerance in bias, while the DWR method requires three levels. Also, when using uniform meshes, more samples on each level are required than in the DWR method in order to reduce the variance.
\begin{table}[h!]
    \centering
    \begin{tabular}{|c|c|c|c|c|}
         \hline
         Level & \# Elems & Cost Per Sample & \# Samples & Var. Per Level \\
         \hline
         0 & 72    & 1     & 1530 & 3.19378E-6 \\
         \hline
         1 &  162  & 3.25  & 29 &  7.93078E-7\\
         \hline
         2 & 288   & 6.25   & 34  & 1.96356E-6 \\
         \hline
          3 & 578  &  12.027  & 10  &  3.15403E-7 \\
         \hline
          4 & 1152 & 24.027   & 10  & 2.41456E-6 \\
         \hline
    \end{tabular}
    \caption{Results for each level of the estimator in example from \S \ref{sec:adv_diff} using uniform meshes with tolerance $\epsilon=10^{-5}$.}
    \label{tab:adv_unif}
\end{table}
\begin{table}[h!]
    \centering
    \begin{tabular}{|c|c|c|c|c|}
    \hline
        Tot. Var. & Squared Bias & MSE & Est. Exp. Val & Tot. Cost \\
    \hline
   8.68040E-6 & 5.87876E-7 & 9.26828E-6 & -0.23215 & 2197.2 \\
    \hline
    \end{tabular}
    \caption{Results of the MLMC estimator in example from \S \ref{sec:adv_diff} using uniform meshes with tolerance $\epsilon=10^{-5}$. }
    \label{tab:adv_unif1}
\end{table}
\begin{table}[h!]
    \centering
    \begin{tabular}{|c|c|c|c|c|}
         \hline
         Level & \# Elems & Cost Per Sample & \# Samples & Var. Per Level \\
         \hline
         0 & 72   & 1     & 979 & 4.98828E-6  \\
         \hline
         1 & 135  & 2.680  & 25 & 8.14877E-8 \\
         \hline
         2 & 291  & 4.972   & 10  & 1.68624E-7 \\
         \hline
    \end{tabular}
    \caption{Results for each level of the estimator in example from \S \ref{sec:adv_diff} with tolerance $\epsilon=10^{-5}$. New meshes are obtained using DWR refinement where the 25\% largest contributions to the error are refined.}
    \label{tab:adv_adapt}
\end{table}
\begin{table}[h!]
    \centering
    \begin{tabular}{|c|c|c|c|c|}
    \hline
        Tot. Var. & Squared Bias & MSE & Est. Exp. Val & Tot. Cost \\
    \hline
     5.23839E-6 & 4.10668E-7 & 5.64906E-6 & -0.22919 & 1120.5 \\
    \hline
    \end{tabular}
    \caption{Results of the MLMC estimator in example from \S \ref{sec:adv_diff} with tolerance $\epsilon=10^{-5}$. New meshes are obtained using DWR refinement where the 25\% largest contributions to the error are refined.}
    \label{tab:adv_adapt1}
\end{table}
\section{Conclusion}
We develop an MLMC algorithm using adjoint-based error analysis to obtain an accurate approximation of the bias and adaptively refine meshes when creating new levels. It is shown that using an adaptive refinement method, either meso-scale or refining regions of largest error, leads to a more cost-effective method than uniform refinement. The advantages of adaptive refinement become more prominent in higher-dimensional problems and in problems where error accumulation is localized. Further improvements can be made to the meso-scale refinement method to better incorporate data from all samples. Other future work will also focus on parabolic PDEs.

\bibliographystyle{siamplain}
\bibliography{refs}

\end{document}